\documentclass{amsart}

\usepackage{enumerate}
\usepackage{amsmath}
\usepackage{amssymb}
\usepackage{amsthm}
\usepackage{mathtools}

\usepackage{graphicx}
\usepackage{type1cm}
\usepackage{eso-pic}
\usepackage{color}
\makeatletter

\usepackage{math}
\usepackage{thm}
\usepackage{ref}

\newcommand{\R}{\ensuremath{\mathbb{R}}}

\newcommand{\D}{\ensuremath{\mathcal{D}}}
\newcommand{\N}{\ensuremath{\mathbb{N}}}

\newcommand{\Div}{\operatorname{div}}
\newcommand{\DIV}{\operatorname{Div}}

\newcommand{\E}{\mathbb{E}}
\newcommand{\F}{\mathbb{F}}
\newcommand{\G}{\mathbb{G}}
\newcommand{\U}{\mathbb{U}}

\newcommand{\pa}{\partial}
\newcommand{\grad}{\operatorname{grad}}

\let\altphi\phi
\renewcommand{\phi}{\varphi}

\title
	[Local Well-Posedness for Relaxational Fluid Vesicle Dynamics]
	{Local Well-Posedness for \\ Relaxational Fluid Vesicle Dynamics}

\author
	[Matthias K{\"o}hne]
	{Matthias K{\"o}hne}

\address
	{Mathematisches Institut,
	 Heinrich-Heine-Universit{\"a}t D{\"u}sseldorf \newline\indent
	 Universit{\"a}tsstr.~1, 40225 D{\"u}sseldorf, Germany}
	
\email
	{koehne@math.uni-duesseldorf.de}

\author
	[Daniel Lengeler]
	{Daniel Lengeler}

\address
	{Fakult{\"a}t f{\"u}r Mathematik,
	 Universit{\"a}t Regensburg \newline\indent
	 Universit{\"a}tsstr.~31, 93053 Regensburg, Germany}
	
\email
	{daniel.lengeler@mathematik.uni-regensburg.de}


\thanks{The second author gratefully acknowledges support by DFG SPP 1506 "Transport Processes at Fluidic Interfaces".}		 


\date
	{\today}

\begin{document}
\begin{abstract}
We prove the local well-posedness of a basic model for relaxational fluid vesicle dynamics by a contraction mapping argument. Our approach is based on the maximal $L_p$-regularity of the model's linearization.
\end{abstract}
\maketitle

\section*{Introduction}
Most biological membranes are composed of a two-layered sheet of phospholipid molecules, a \emph{lipid bilayer}, which is immersed in water. Due to hydrophobic effects, these membranes tend to avoid open edges and form closed configurations called \emph{vesicles}. Since the ratio of membrane thickness to vesicle  diameter is very small, typically of the order $10^{-4}$, vesicles can be described as two-dimensional surfaces embedded in three-dimensional space. Due to osmotic effects and a very low solubility of the phospholipids, the area and enclosed volume of such a vesicle are practically fixed. Hence, vesicle configurations 
are not determined by a surface tension but rather by a bending elasticity. A basic model for such an elastic energy is given by the \emph{Canham-Helfrich energy}
\begin{equation*}
 F(\Gamma)=\frac\kappa 2\int_\Gamma (H-C_0)^2\ dA;
\end{equation*}
see \cite{canham,helfrich,evans}. Here, $\Gamma$ is the two-dimensional, closed surface representing the membrane, $H$ denotes twice its mean curvature, $\kappa$ is the bending rigidity, and $C_0$ is the \emph{spontaneous curvature}, which is supposed to reflect a chemical asymmetry of the membrane or its environment; both $\kappa$ and $C_0$ are assumed to be constant in the following. Usually the lipid bilayer is in a fluid state, allowing the monolayers to freely flow laterally and to slip over each other while the membrane retains its transverse structure. In our basic model we take into account this fluidity while neglecting the bilayer architecture of the membrane. More precisely, we study a single homogeneous Newtonian surface fluid (see \cite{scriven,aris}) subject to additional stresses induced by the Canham-Helfrich energy and interacting with a homogeneous Newtonian bulk fluid. The full system reads as follows:
\begin{equation}\label{eqn:nonrex}
 \begin{aligned}
 \rho_b\frac{Du}{Dt}&=\Div S&&\mbox{ in }\Omega\setminus\Gamma_t,\\
 \Div u&=0&&\mbox{ in }\Omega\setminus\Gamma_t,\\
 V_t = u \cdot \nu_t, \quad [\![u]\!] = 0, \quad \rho\frac{Du}{Dt}&=\DIV T + [\![S]\!]\nu_t &&\mbox{ on }\Gamma_t,\\
 \DIV u&=0&&\mbox{ on }\Gamma_t,\\
 u&=0&&\mbox{ on }\pa\Omega,\\
 u(0)&=u_0&&\mbox{ in }\Omega\setminus\Gamma_0,\\
 \Gamma_t|_{t=0}&=\Gamma_0.
 \end{aligned}
 \end{equation}
Here, $\Omega$ is a smooth bounded domain in $\R^3$ containing a homogeneous Newtonian fluid and a closed moving vesicle $\Gamma_t$, $\nu_t$ is the outer unit normal on $\Gamma_t$, $u$ is the velocity of the bulk fluid in $\Omega\setminus\Gamma_t$ and the velocity of the surface fluid on $\Gamma_t$ which are assumed to coincide on $\Gamma_t$, $\rho_b$ and $\rho$ denote the bulk and the surface mass density, respectively, $Du/Dt$ is the fluid particle acceleration, $S=2\mu_b Du - \pi I$ is the Newtonian bulk stress tensor with the constant dynamic viscosity $\mu_b$ of the bulk fluid, the symmetric part $Du$ of the gradient of $u$, and the bulk pressure $\pi$, 
$V_t$ is the speed of normal displacement of $\Gamma_t$,
$[\![u]\!]$ and $[\![S]\!]$ are the the jump of the velocity and the bulk stress tensor, respectively, across the membrane (subtracting the \emph{outer limit} from the \emph{inner limit}), $\DIV$ is the surface divergence (see below), and $T={^fT}+{^eT}$ is the surface stress tensor which is composed of a fluid part ${^fT}$ and an elastic part ${^eT}$ induced by the Canham-Helfrich energy. More precisely, in coordinates we have $^fT^i_\alpha={^f\tilde T_\alpha^\beta}\pa_\beta^i$ with (cf. \cite{scriven,aris,lengeler})
\[^f\tilde T_\alpha^\beta=-q\,\delta_\alpha^\beta+2\mu\,(\D u)_\alpha^\beta=-q\,\delta_\alpha^\beta+\mu\,g^{\beta\gamma}(v_{\alpha;\gamma}+v_{\gamma;\alpha}-2w\,k_{\alpha\gamma})\]
and 
\[^{e}T^i_\alpha=\kappa\,\big((H-C_0)^2/2\,\pa_\alpha^i-(H-C_0)\,k_\alpha^\beta\pa_\beta^i-(H-C_0)_{,\alpha}\nu^i_t\big).\]
Here, $q$ is the surface pressure acting as a Lagrange multiplier with respect to the constraint $\DIV u=0$, $\mu$ is the constant dynamic viscosity of the surface fluid, $\D u$ is the \emph{surface rate-of-strain tensor}, $k$ is the second fundamental form of $\Gamma_t$, $\partial_\alpha$ denotes the $\alpha$-th coordinate vector field, and the semicolon denotes covariant differentiation while the comma indicates usual partial differentiation. Furthermore, on $\Gamma_t$ we decomposed the function $u=v+w\,\nu_t$ into its tangential and its normal part. Throughout the paper, latin indices refer to Cartesian coordinates in $\R^3$ while greek indices refer to arbitrary coordinates on $\Gamma_t$. In particular, we note that the surface stress tensors are instances of \emph{hybrid tensor fields} (see \cite{scriven,aris}) taking a tangential direction and returning a force density that is, in general, not tangential. The surface divergences for the non-tangential vector field $u$ and the hybrid tensor field $T$ 
can 
be written as
\begin{equation*}
\begin{aligned}
\DIV u &=g^{\alpha\beta}\langle\pa_\alpha u,\pa_\beta\rangle_{e},\\
(\DIV T)^i&=g^{\alpha\beta} T^i_{\alpha;\beta},
\end{aligned}
\end{equation*}
where $g$ denotes the Riemannian metric on $\Gamma_t$ induced by the Euclidean metric $e$ in $\R^3$, and the semicolon denotes the corresponding covariant differentiation of the covectors $(T^i_{\alpha})_{\alpha=1,2}$ (for fixed $i$). The computations in \cite{lengeler} showed that
\begin{equation}\label{fullt}
\begin{aligned}
\DIV u=&\Div_g v - w\, H,\\
 \DIV T = &-\grad_g q -q\,H\nu_t + \mu\,\big(\Delta_g v + \grad_g(w\,H) + K v -2\Div_g(w\,k)\big)\\ &  +2\mu\big(\langle\nabla^g v,k\rangle_g-w\,(H^2-2K)\big)\nu_t\\
& -\kappa\big(\Delta_g H + H(H^2/2-2K)+C_0(2K-HC_0/2)\big)\nu_t.
\end{aligned}
\end{equation}
Here, $K$ is the Gauss curvature, $\grad_g$, $\Div_g$, $\nabla^g$, $\Delta_g$ denote the differential operators (acting on tangential tensor fields) corresponding to the metric $g$, and, with a slight abuse of notation, we write $\langle\nabla^g v,k\rangle_g$ for the contraction of the tensor fields $\nabla^g v$ and $k$ using $g$. Furthermore, we saw in \cite{lengeler} that both the bulk and the surface Reynolds number usually are very small, typically of the order $10^{-3}$. Hence, neglecting the inertial terms in \eqref{eqn:nonrex} we arrive at the following set of equations describing purely relaxational fluid vesicle dynamics:
\begin{equation}\label{eqn:final}
\begin{aligned}
\Div S&=0 &&\mbox{ in }\Omega\setminus \Gamma_t,\\
\Div u&=0&&\mbox{ in }\Omega\setminus \Gamma_t,\\
V_t = u \cdot \nu_t, \quad [\![u]\!] = 0, \quad \DIV {^fT}+[\![S]\!]\nu_t&=-\DIV {^eT}&&\mbox{ on }\Gamma_t,\\
\DIV u&=0&&\mbox{ on }\Gamma_t,\\
u&=0&&\mbox{ on }\pa\Omega,\\
\Gamma_t|_{t=0}&=\Gamma_0.
\end{aligned}
\end{equation}
At first sight, one might think that there is no dynamical component left in the system. However, this is not the case. Note that $\DIV {^eT}$ can be computed from $\Gamma_t$ alone. Hence, we have to solve the Stokes-type system defined by the left hand side of \eqref{eqn:final} with $-\DIV {^eT}$ as a right hand side for the fluid velocity $u$. Then, the normal part $w$ of $u$ on $\Gamma_t$ tells us how the vesicle will move in the next instant. It is easy to conclude from \eqref{eqn:final}$_{2,4}$ that the area and the enclosed volume of each connected component $\Gamma^i_t$ of $\Gamma_t$ are preserved under this flow; see \cite{lengeler}. Hence, the phase space $N$ of \eqref{eqn:final} consists of the embedded surfaces $\Gamma\subset\Omega$ of fixed area and enclosed volume. As \eqref{fullt}$_2$ indicates (see also \cite{lengeler}) we have $-\DIV {^eT}=\grad_{L_2}F_{\Gamma_t}\,\nu_t$, where $\grad_{L_2}F_{\Gamma_t}$ denotes $L_2$-gradient of the Canham-Helfrich energy at the point $\Gamma_t$. Hence, 
we 
note that compared to the classical Canham-Helfrich flow, that is, the $L_2$-gradient flow of the 
Canham-Helfrich energy with prescribed enclosed volume and area, there is an additional Neumann-to-Dirichlet-type operator involved here, mapping $\grad_{L_2}F_{\Gamma_t}$ to $w$; since $\grad_{L_2}F_{\Gamma_t}$ is a fourth order operator in $\Gamma_t$, the mapping $\Gamma_t\mapsto w$ can be considered as a nonlinear, nonlocal pseudo-differential operator of third order. Furthermore, we saw in \cite{lengeler} that \eqref{eqn:final} can be considered as a gradient flow with respect to a suitable Riemannian metric on $N$, leading in particular to the energy identity
\begin{equation}\label{enid}
\frac{d}{dt} F(\Gamma_t)=-2\mu_b\int_{\Omega\setminus \Gamma_t}|Du|^2_{e}\,dx-2\mu\int_{\Gamma_t}|\D u|^2_g\,dA.
\end{equation}
We will not make use of the gradient flow structure in the present article. However, it turns out to be useful for the proof of asymptotic stability of local minimizers of the Canham-Helfrich energy; this is done in \cite{lengelerstability} by using a \L ojasiewicz-Simon inequality. Finally, we showed in \cite{lengeler} that the equilibria $\Gamma$ of \eqref{eqn:final} satisfy
\begin{equation*}
 \grad_{L^2}F_{\Gamma}+[\![\pi]\!]+q\,H=0.
\end{equation*}
This is the \emph{Helfrich equation} with the pressure jump and the surface pressure acting as Lagrange multipliers with respect to the volume and area constraints.

Not much rigorous analysis has been done on the dynamics of fluid vesicles. Concerning the \emph{Canham-Helfrich flow}, a partial local well-posedness result has been shown in \cite{nagasawa12}. There exist further results \cite{nagasawa06,wheeler1,wheeler2} concerning a Helfrich-type flow where the Lagrange parameters instead of volume and area are prescribed and which consequently should not be related directly to fluid vesicles. In \cite{Zhang} local-in-time existence and uniqueness for a homogeneous Newtonian surface fluid subject to Canham-Helfrich stresses is shown. While the bulk fluid is neglected the authors keep the inertial term in the equations for the surface fluid, yielding a kind of dissipative fourth order wave-type equation. In \cite{Shkoller} local-in-time existence and uniqueness of a homogeneous Newtonian bulk fluid with inertial term interacting with a compressible, inviscid surface fluid without inertial term is shown, the membrane model being rather
non-standard. Since the authors work in the $L_2$-scale they have to deal with solutions of higher regularity, making the analysis rather involved. Furthermore, they work in the 
Lagrangian picture, leading to problems with the tangential degeneracy of the elliptic operator arising from the elastic stresses within the membrane. By working in an $L_p$-scale and using the \emph{Hanzawa transform} instead of the Lagrangian picture, we are able to present a simplified analysis based on the theory of maximal $L_p$-regularity along with localization and transformation techniques.

The present article continues our analysis of a basic model of fluid vesicle dynamics that was started in \cite{lengeler}, where a thorough $L_2$-analyis of the Stokes-type system defined by the left hand side of \eqref{eqn:final} is performed. We will make extensive use of these results in the present article. Furthermore, we refer the reader to \cite{lengeler} and the references therein for a detailed introduction to the physics and mathematics of fluid vesicles; in particular, we refer to \cite{seifert97} for a rather comprehensive treatment of the physics of equilibrium configurations.

This paper is organized as follows. In \Secref{mainresult} we present our main result. In order to construct local solutions, we shall employ the standard procedure of approximating the nonlinear evolution by some appropriate linear evolution. To this end, in \Secref{hanzawa}, we transform our system to a fixed domain using the Hanzawa transform and extract the linearization of the resulting system. In \Secref{linear} we prove that the linearization has the property of maximal $L_p$-regularity. This is first done for the case of a double half-space, to which, then, the general case is reduced by localization and transformation techniques. Finally, in \Secref{iteration}, we prove our main result, using the contraction mapping principle.

Before we proceed, let us fix some notation. Throughout the article, let $\Omega\subset\R^3$  be a smooth bounded domain and $\Gamma\subset\Omega$ a smooth, closed surface with outer unit normal $\nu$. We write $\Gamma^i$, $i=1,\ldots,m$, for the connected components of $\Gamma$, $\Omega^i$ for the open set enclosed by $\Gamma^i$, and let \[\Omega^0:=\Omega\setminus\big(\bigcup_{i=1}^m \Gamma^i\cup\Omega^i\big).\]
We denote by $P_\Gamma$ the field of orthogonal projections onto the tangent spaces of $\Gamma$ while $[u]_\Gamma$ denotes the trace of the bulk field $u$ on $\Gamma$; however, when there is no danger of confusion we will sometimes omit the brackets. Furthermore, (apart from \Appref{app:cov}) we write $e$ for the Euclidean metric in $\R^3$ and  $g$ for the metric on $\Gamma$ induced by $e$. We also use the notation $u\cdot v$ instead of $\langle u,v\rangle_e$ for $u,v\in\R^3$. Moreover, we write $k$, $H$, and $K$ for the second fundamental form, twice the mean curvature, and the Gauss curvature of $\Gamma$ with respect to $e$, respectively. With a slight abuse of notation we use same symbol $k$ also to denote the Weingarten map, that is, in coordinates we write $k_{\alpha\beta}$ and $k_\alpha^\beta$. Furthermore, for any metric $\tilde e$ on an arbitrary manifold, we write $^{\tilde e}\Gamma_{ij}^k$, 
$\nabla^{\tilde e}$, $\grad_{\tilde e}$, $K_{\tilde e}$, etc. 
for the associated Christoffel symbols,
 differential operators, and curvature terms, and we use the abbreviations $\Gamma_{ij}^k:={^e\Gamma}_{ij}^k$, $\nabla:=\nabla^e$, $\grad:=\grad_{e}$, etc. for the corresponding Euclidean objects. When working in coordinates and 
confusion about the underlying metric can be ruled out, we use the semicolon to separate the indices coming from covariant differentiation from the original indices; for instance, for a covector field $\omega$ we write $(\nabla^{\tilde e} \omega)_{ij}=\omega_{i;j}$. We denote by $r(a)$ generic tensor fields that are polynomial or analytic functions of their argument $a$. Furthermore, for tensor fields $r_1$ and $r_2$ we write $r_1*r_2$ for any tensor field that depends in a bilinear way on $r_1$ and $r_2$, and we use the abbreviations $r*(r_1,\ldots,r_k)=r*r_1 + \ldots + r*r_k$ and $r^k=r*\ldots*r$ (with $k$ factors on the right hand side).
For $p\in (1,\infty)$, $k\in\N$, and $s\in\R_+\setminus\N$ we denote by $H^k_p$ the usual Sobolev spaces and by $W^s_p$ the Sobolev-Slobodeckij spaces. For an arbitrary smooth, $d$-dimensional Riemannian manifold $(M,\tilde e)$ the norm of the latter spaces is given by
\[\|T\|_{W^s_p(M)}=\|T\|_{H^k_p(M)}+|(\nabla^{\tilde e})^k T|_{W^{s-k}_p(M)},\]
where $k$ is the largest integer smaller than $s$ and 
\[|(\nabla^{\tilde e})^k T|_{W^{s-k}_p(M)}^p:=\int_M\int_M \frac{|(\nabla^{\tilde e})^k T(x)-(\nabla^{\tilde e})^k T(y)|^p_e}{d_{\tilde e}(x,y)^{d+(s-k)p}}\,dV_{\tilde e}(x)\,dV_{\tilde e}(y).\]
In this formula $d_{\tilde e}$ is the Riemannian distance function while $dV_{\tilde e}$ is the volume element corresponding to $\tilde e$. Finally, let the homogeneous spaces $\dot H^k_p(M)$ and $\dot W^s_p(M)$ consist of all locally integrable tensor fields such that $\nabla^k T\in L_p(M)$ and $|\nabla^k T|_{W^{s-k}_p(M)}<\infty$ (where $k$ is the largest integer smaller than $s$), respectively.

\section{Main result}\seclabel{mainresult}
We fix a smooth bounded domain $\Omega\subset\R^3$ and a smooth, closed surface $\Gamma\subset\Omega$. We denote by $S_{\alpha}$, $\alpha>0$, the open set of points in $\Omega$ whose distance from $\Gamma$ is less than $\alpha$. It's a well-known fact from elementary differential geometry that there exists some $\gamma>0$ such that the mapping
\begin{equation*}
 \begin{aligned}
\Lambda: \Gamma\times (-\gamma,\gamma)\rightarrow
S_{\gamma},\
(x,d)\mapsto x + d\,\nu(x)  
 \end{aligned}
\end{equation*}
is a diffeomorphism. For functions $h:\Gamma\rightarrow (-\kappa,\kappa)$ we define
\begin{equation*}
	\Gamma_h := \{\,\Lambda(x,\,h(x))\,|\,x\in\Gamma\,\},
\end{equation*}
and we write $x\mapsto(\tau(x),d(x)): S_\gamma \longrightarrow \Gamma \times (-\gamma,\gamma)$ for the inverse map $\Lambda^{-1}$,
i.\,e.\ we denote by $\tau: S_\gamma \longrightarrow \Gamma$ the metric projection onto $\Gamma$,
and by $d: S_\gamma \longrightarrow (-\gamma,\gamma)$ the signed distance from $\Gamma$,
which are both well-defined within $S_\gamma$ by choice of $\gamma > 0$.
For a time-dependent closed surface $\Gamma_t\subset\Omega$ and time-dependent, integrable, scalar functions $q,\pi$ defined on $\Gamma_t$ and in $\Omega$, respectively, we consider the gauge conditions
\begin{equation}\label{compat1}
 \int_{\Gamma^i_t}q(t,\cdot\,)/H\,dA+\int_{\Omega^i_t}\pi(t,\cdot\,)\,dx = 0\quad\text{for each $\Gamma^i_t$ that is a round sphere},
\end{equation}
where $\Omega_t^i$ denotes the open set enclosed by $\Gamma_t^i$, and
\begin{equation}\label{compat2}
 \int_{\Omega}\pi(t,\cdot\,)\,dx=0.
\end{equation}
Note that condition \eqref{compat1} is a consequence of the divergence constraint on $\Gamma^i_t$ and in $\Omega \setminus \Gamma_t$, provided that $\Gamma^i_t$ is a CMC surface, i.\,e.\ provided $\Gamma^i_t$ is a round sphere.
Now, we are ready to state our main result. Let $\mu_b,\mu>0$.
\begin{theorem}\label{thm}
Assume that $\Gamma$ contains no round sphere, and let $p\in (3,\infty)\setminus\{4\}$. For sufficiently small $\epsilon>0$ there exists a time $T>0$ such that for all height functions $h_0\in \bar B_\epsilon(0)\subset W^{5-4/p}_p(\Gamma)$ there exists a solution of \eqref{eqn:final} in the time interval $I:=(0,T)$ with initial value $\Gamma_{h_0}$. This solution is given by $\Gamma_t=\Gamma_{h(t)}$ for a height function 
 \begin{equation*}
 h\in H^1_p(I,\,W^{2 - 1/p}_p(\Gamma)) \cap L_p(I,\,W^{5 - 1/p}_p(\Gamma))
\end{equation*}
such that $\|h\|_{L_\infty(I\times \Gamma)}<\gamma$ and by measurable hydrodynamic fields $u$, $\pi$ defined in $\Omega\setminus\Gamma_t$ and $q$ defined on $\Gamma_t$ for almost all $t\in I$ such that the functions
 \begin{equation*}
 \|u(t,\cdot\,)\|_{H^{2}_p(\Omega\setminus\Gamma_t)}^p,\ \|P_{\Gamma_t} [u(t,\cdot\,)]_{\Gamma_t}\|_{H^{2}_p(\Gamma_t)}^p,\ \|\pi(t,\cdot\,)\|_{H^{1}_p(\Omega\setminus\Gamma_t)}^p,\
 \|q(t,\cdot\,)\|_{H^{1}_p(\Gamma_t)}^p
\end{equation*}
are integrable in $I$,
and such that \eqref{compat1} and \eqref{compat2} hold for almost all $t$; the solution is unique in this class. Moreover, the map
\begin{equation*}
\bar B_\epsilon(0)\subset W^{5-4/p}_p(\Gamma)\rightarrow H^1_p(J,\,W^{2 - 1/p}_p(\Gamma)) \cap L_p(J,\,W^{5 - 1/p}_p(\Gamma)),\quad h_0\mapsto h
\end{equation*}
is Lipschitz continuous.

Finally, if $\Gamma$ consists only of round spheres, then a solution of \eqref{eqn:final} is given by the constant-in-time solution with $\Gamma_t=\Gamma$, $u=0$, and suitably chosen pressure functions $\pi$ and $q$ which are constant in each connected component of $\Omega$ and $\Gamma$, respectively; this solution is unique in the class given in the first part of the theorem. In particular, the problem is globally well-posed in this trivial case.
\end{theorem}

In general, when dealing with the continuous dependence part of (local) well-posedness the question arises which perturbations should be included in the analysis. For macroscopic physical systems it seems reasonable to consider those perturbations which are accessible by thermal fluctuations; thus, in our case area and enclosed volume of each connected component $\Gamma^i_{h_0}$, $i=1,\ldots,m$, of $\Gamma_{h_0}$ should be conserved. Concerning the second part of the above theorem, note that consequently the only admissible perturbations of a round sphere are translations of this sphere. On the other hand, the first part of the above theorem is slightly more general in that it deals with a larger class of perturbations.
Note that the tangential part of the bulk velocity trace on $\Gamma_t$ exhibits an increased spatial regularity, which is to be expected in view of the appearance of the Laplace-Beltrami operator in the transformed equation $(\ref{eqn:tildesystem})_3$; cf.~also the symbolic analysis on page~\pageref{eqn:tpzv}.
Also note that the case $p = 4$ is excluded for notational convenience, since in this case Besov spaces would have to be introduced for the initial data.

So far, we cannot prove (local) well-posedness of our system in the case that $\Gamma$ contains both round spheres and connected components that are not round spheres. The reason for this is a technicality in the iterative construction process of the solutions which is related to the different degrees of gauge freedom for round spheres on the one hand and other configurations on the other hand; see, in particular, the remark following the proof of Theorem \ref{thm}. 

The conditions \eqref{compat1} and \eqref{compat2} on $\pi$ and $q$ provide a gauge fixing; as is typical for Stokes-type equations the pressure functions in \eqref{eqn:final} are not uniquely determined.
\begin{definition}\label{def}
For fixed $t \in \bar I$ we define the space $U_p(\Gamma_t) \subset H^1_p(\Omega\setminus\Gamma_t) \times H^1_p(\Gamma_t)$ as follows: $(\pi,q) \in U_p(\Gamma_t)$, if and only if for all $i=1,\ldots,m$ we have
\begin{itemize}
 \item[(i)] $\pi=\kappa_i$ in $\Omega^i$, $\pi=\kappa_0$ in $\Omega^0$, $q=\kappa^i$ on $\Gamma^i$ with $\kappa_i,\kappa_0,\kappa^i\in\R$
 \item[(ii)] If $\Gamma^i_t$ is a round sphere with $H$ denoting twice the mean curvature, then $\kappa_i-\kappa_0=\kappa^i H$.
 \item[(iii)] If $\Gamma^i_t$ is not a round sphere, then $\kappa^i=0$ and $\kappa_i=\kappa_0$.
\end{itemize}
\end{definition}
It is not hard to see that 
\begin{equation*}
 \begin{aligned}
\big(H^1_p(\Omega\setminus\Gamma_t)\times H^1_p(\Gamma_t)\big)/U_p(\Gamma_t)
\simeq \big\{(\pi,q)\in H^1_p(\Omega\setminus\Gamma_t)\times H^1_p(\Gamma_t)\,|\,\mbox{\eqref{compat1},\eqref{compat2} hold}\big\} 
\end{aligned}
\end{equation*}
cf.~Section 3.1 in \cite{lengeler}. Hence, the subspace $U_p(\Gamma_t)$ characterizes the gauge freedom of the pressure functions.
Concerning the gauge fixing condition \eqref{compat1} we note that a connected component $\Gamma_t^i$, $i=1,\ldots,m$, of $\Gamma_t$ is a round sphere for some $t\in \bar I$ if and only if this is the case for all $t\in \bar I$ since its area and enclosed volume are fixed. 

Assuming the reference surface (respectively, initial surface in the case $h_0=0$) $\Gamma$ to be of class $W^{6-1/p}_p$ would be sufficient as a detailed analysis of the nonlinearities appearing in \Secref{hanzawa} shows. By Theorem 4.10.2 in \cite{amann95} and the theorem in Section 7.4.4 of \cite{Triebel:Function-Spaces-2} (note that $W^s_p(\Gamma)=B^s_{pp}(\Gamma)$ for non-integer $s$; cf. \cite{lengeler}) the time trace
\[H^1_p(I,\,W^{2 - 1/p}_p(\Gamma)) \cap L_p(I,\,W^{5 - 1/p}_p(\Gamma))\rightarrow W^{5-4/p}_p(\Gamma),\quad h\mapsto h(0)\] is well-defined and surjective. In this sense the regularity of $h_0$ in the preceding theorem is optimal. However, so far, we are not able to prove the well-posedness for arbitrary (apart from spheres) initial surfaces $\Gamma$ of class $W^{5-4/p}_p$. The canonical way to do so is to approximate such an initial surface sufficiently well by some smooth surface and then apply the preceding theorem. However, it seems that the $\epsilon$ in the assertion (being related to the norm of the solution operator of the linearization with respect to the reference surface) does not depend in a continuous way on the reference surface in the $W^{5-4/p}_p$-topology. However, it should be possible to lower the regularity assumption on $\Gamma$ below $W^{6-1/p}_p$ by proving such a continuity result in a sufficiently strong topology.

Finally, we note that a similar result as Theorem \ref{thm} result can be shown for $\mu=0$; in this case, of course, the additional regularity of the tangential velocity on the membrane is not present.


\section{Linearization}
\seclabel{hanzawa}
In this section we employ the classical Hanzawa transform to map the time-dependent domains $\Gamma_t$ and $\Omega\setminus \Gamma_t$ to the fixed domains $\Gamma$ and $\Omega\setminus \Gamma$, respectively. Using this diffeomorphism we translate the system \eqref{eqn:final} into a quasi-linear system on fixed domains and then extract its linearization. It is crucial, however, to transform the equations in a geometrically consistent way, namely, to take the geometric pull-back of the fields involved. This ensures that the tangential part of the velocity field on $\Gamma_t$, which is smoothed by membrane viscosity, remains tangential, and thus is smoothed in the linearization, too.\footnote{At first sight, one might want to transform the equation in such a way that the normal part of velocity on $\Gamma_t$ remains normal. However, the construction of a suitable modification of the classical Hanzawa transform turns out to be rather technical (see \cite{meyer} for an elegant method), 
and, in our case, the need for this property can be avoided by relaxing the relation $u=v+w\,\nu_t$ in the linearization; see below.} 

Recall the notation from the beginning of \Secref{mainresult}. For sufficiently regular $h:\Gamma\rightarrow (-\gamma,\gamma)$, we choose the real-valued function $\beta\in C^\infty(\R)$ to be $0$ in neighbourhoods of $-1$ and $1$, and $1$ in a neighborhood of $0$, and assume that $|\beta'|< \gamma/\|h\|_{L^\infty(\Gamma)}$ on $\Gamma$. Then, the Hanzawa transform $\Phi_h:\overline\Omega\rightarrow\overline\Omega$ is defined in the following way: While, in $\overline\Omega\setminus S_\gamma$, we let $\Phi_h$ be the identity, we define $\Phi_h$ in $S_\gamma$ by
\begin{equation*}
 \begin{aligned}
x&\mapsto x + \nu(\tau(x))\,h(\tau(x))\,\beta(d(x)/\gamma).
\end{aligned}
\end{equation*}
Then we have $\Phi_h(\Gamma)=\Gamma_h$, and it is not hard to see that $\Phi_h:\overline\Omega\rightarrow\overline\Omega$ and $\varphi_h:=\Phi_h|_{\Gamma}:\Gamma\rightarrow\Gamma_h$ are diffeomorphisms; see, for instance, \cite{lengelerphd}. For a given time-dependent height function $h$, we write $\Phi_t:=\Phi_{h(t)}$, $\varphi_t:=\varphi_{h(t)}$, and $\Gamma_t:=\Gamma_{h(t)}$.

Separating the tangential and the normal part of \eqref{eqn:final}$_3$, we obtain
\begin{equation}\label{eqn:divt}
\begin{aligned}
&-\grad_g q + \mu\,\big(\Delta_g v + \grad_g(w\,H) + K v -2\Div_g(w\,k)\big)+2\mu_b[\![Du]\!]\nu=0,\\
&-q\,H + 2\mu\big(\langle\nabla^g v,k\rangle_g-w\,(H^2-2K)\big)-[\![\pi]\!]\\
&\quad\quad\quad\quad\quad\quad\quad=\kappa\big(\Delta_g H + H(H^2/2-2K)+C_0(2K-HC_0/2)\big).
\end{aligned}
\end{equation}
Note that $[\![Du]\!]\nu$ is tangential due to the incompressibility constraint. Indeed, for any vector $X$ on $\Gamma$ we have $[\![(X\cdot\nabla)u]\!]\cdot\nu=0$. If $X$ is tangential, we even have $[\![(X\cdot\nabla)u]\!]=0$ since $u$ is continuous across $\Gamma$. But then, choosing an orthonormal basis $\nu,e_1,e_2$ at some arbitrary point on $\Gamma$, from $\Div u=0$ we deduce that
\[[\![(\nu\cdot\nabla)u]\!]\cdot\nu=-[\![(e_1\cdot\nabla)u]\!]\cdot e_1-[\![(e_2\cdot\nabla)u]\!]\cdot e_2=0.\]

Let us now transform the system \eqref{eqn:final} to the fixed domains $\Omega\setminus \Gamma$ and $\Gamma$ and then extract its linearization. We minimize the computations by exploiting the fact that the system \eqref{eqn:final} on the time-dependent domains is equivalent to a system on the fixed domains with a time-dependent Riemannian metric. The diffeomorphism $\Phi_t$ induces the time-dependent metric $\tilde e=\tilde e_t:=\Phi^*_t e$ on $\Omega$. We denote the restriction of $\tilde e$ to $\Gamma$ by $\tilde g$. Note that $\Phi_t:(\Omega,\tilde e_t)\rightarrow (\Omega,e)$ and $\varphi_t:(\Gamma,\tilde g_t)\rightarrow (\Gamma_t,g)$ are isometries. Let us denote the pullbacks of the involved fields by $\tilde u:=\Phi_t^*u$, $\tilde\pi:=\Phi_t^*\pi$, $\tilde v:=\Phi_t^*v$, $\tilde w:=\Phi_t^*w$, and $\tilde q:=\Phi_t^*q$. By exploiting naturality of covariant differentiation under isometries (cf. \cite{lengeler}), from \eqref{eqn:final}, \eqref{eqn:divt}, $u=v+w\,\nu_t$ on $\Gamma_t$, and $\pa_t h\,\nu=(
w\,\nu_t)\circ \varphi_t$ we obtain
\begin{equation}\label{eqn:tildesystem}
\begin{aligned}
\mu_b\Delta_{\tilde e}\tilde u-\grad_{\tilde e}\tilde\pi&=0&&\mbox{ in }\Omega\setminus\Gamma,\\
\Div_{\tilde e} \tilde u&=0&&\mbox{ in }\Omega\setminus\Gamma,\\
\mu\big(\Delta_{\tilde g} \tilde v + \grad_{\tilde g}(\tilde w\,H_{\tilde e}) + K_{\tilde g}\,\tilde v -2\Div_{\tilde g}(\tilde w\,k_{\tilde e})\big)&\\
-\grad_{\tilde g} \tilde q  
+2\mu_b[\![D^{\tilde e}\tilde u]\!]\nu_{\tilde e}&=0&&\mbox{ on } \Gamma,\\
2\mu\big(\langle\nabla^{\tilde g} \tilde v,k_{\tilde e}\rangle_{\tilde g}
-\tilde w\,(H^2_{\tilde e}-2K_{\tilde g})\big)-\tilde q\,H_{\tilde e}-[\![\tilde\pi]\!]& &&\\
-\kappa\big(\Delta_{\tilde g} H_{\tilde e} + H_{\tilde e}(H_{\tilde e}^2/2-2K_{\tilde g})+C_0(2K_{\tilde g}-H_{\tilde e}C_0/2)\big)&=0&&\mbox{ on }\Gamma,\\
\Div_{\tilde g} \tilde v-\tilde w\, H_{\tilde e}&=0&&\mbox{ on }\Gamma,\\
\tilde u-\tilde v-\tilde w\,\nu_{\tilde e}&=0&&\mbox{ on }\Gamma,\\
\partial_t h\,\langle\nu,\nu_{\tilde e}\rangle_{\tilde e}-\tilde w&=0&&\mbox{ on }\Gamma.
\end{aligned}
\end{equation}
Here, most of the geometric quantities have to be taken with respect to the perturbed metrics and, hence, are indexed by $\tilde e$ and $\tilde g$, respectively. Now, we take the point of view of $\tilde e$ being a (small) perturbation of $e$. The results from \Appref{app:cov} show (cf. the proof of Theorem 3.13  in \cite{lengeler}) that \eqref{eqn:tildesystem} can be written in the form
\begin{equation}\label{eqn:tildesystem2}
\begin{aligned}
\mu_b\Delta\tilde u-\grad\tilde\pi&=N_1&&\mbox{ in }\Omega\setminus\Gamma,\\
\Div \tilde u&=N_2&&\mbox{ in }\Omega\setminus\Gamma,\\
\mu\big(\Delta_{g} \tilde v + \grad_{g}(\tilde w\,H) + K\,\tilde v -2\Div_{g}(\tilde w\,k)\big)&\\
-\grad_{g} \tilde q  
+2\mu_b[\![D\tilde u]\!]\nu&=N_3^\top&&\mbox{ on } \Gamma,\\
2\mu\big(\langle\nabla^{g} \tilde v,k\rangle_{g}
-\tilde w\,(H^2-2K)\big)-\tilde q\,H-[\![\tilde\pi]\!]-Ah&=N_3^\perp&&\mbox{ on }\Gamma,\\
\Div_{g} \tilde v-\tilde w\, H&=N_4&&\mbox{ on }\Gamma,\\
\tilde u-\tilde v-\tilde w\,\nu&=N_5&&\mbox{ on }\Gamma,\\
\partial_t h-\tilde w&=N_6&&\mbox{ on }\Gamma
\end{aligned}
\end{equation}
with
\begin{equation*}
\begin{aligned}
N_1&=(\tilde e-e)*(\mu_b\nabla^2 \tilde u,\grad\tilde\pi) + \mu_b\, r(\tilde e)*\big((\nabla^2\tilde e,(\nabla\tilde e)^2)*\tilde u+\nabla\tilde e*\nabla \tilde u\big),\\
N_2&=r(\tilde e)*\nabla\tilde e*\tilde u,\\
N_3^\top&=(\tilde e-e)*(\mu(\nabla^g)^2 \tilde v,\grad_g \tilde q) + \mu_b\,r(\tilde e)*\big((\tilde e-e)*[\nabla \tilde u]+\nabla\tilde e*[\tilde u]\big)\\
&\quad+ \mu\,r(\tilde e)*\big((\tilde e-e)*k^2,(\tilde e-e)*\nabla k,k*\nabla\tilde e,\nabla^2\tilde e,(\nabla\tilde e)^2\big)*[\tilde u]\\
&\quad + \mu\,r(\tilde e)*\big((\tilde e-e)*k,\nabla\tilde e\big)*[\nabla\tilde u]\\
N_3^\perp&=r(\tilde e)*\big((\tilde e-e)*k,\nabla\tilde e\big)\,\tilde q+\mu\,r(\tilde e)*\big((\tilde e-e)*k,\nabla\tilde e\big)*\nabla^g\tilde v \\
&\quad + \mu\,r(\tilde e)*\big((\tilde e-e)*k^2,k*\nabla\tilde e,(\nabla\tilde e)^2\big)*[\tilde u]\\
&\quad+\kappa\big(\Delta_g H+H(H^2/2-2K)+C_0(2K-HC_0/2)\big) + \kappa Q(h),\\
N_4&=r(\tilde e)*\big((\tilde e-e)*k,\nabla\tilde e\big)*[\tilde u],\\
N_5&=r(\tilde e)*(\tilde e- e)\,\tilde w,\\
N_6&=r(\tilde e)*(\tilde e-e)\,\partial_t h,
\end{aligned}
\end{equation*}
where $Ah=\kappa\big(\Delta_g^2 h + (a^{\alpha\beta}h_{,\alpha})_{;\beta} + b\,h\big)$ is the linearization at $h\equiv 0$ of $\grad_{L_2}F_{\Gamma_h}$ with 
\begin{equation*}
\begin{aligned}
a^{\alpha\beta}&=(H^2/2-4K+2HC_0-C_0^2/2)g^{\alpha\beta}+2(H-C_0)k^{\alpha\beta},\\
b&=2k^{\alpha\beta}H_{;\alpha\beta} + \Delta_g(H^2-2K)
+ H_{,\alpha}H_,^{\,\alpha}+3H^4/2 - 7KH^2\\
&\quad+4K^2 +2C_0KH-C_0^2/2\,H^2+C_0^2K,
\end{aligned}
\end{equation*}
see \cite{lengelerstability}, and $\kappa Q(h)=\phi_t^*(\grad_{L_2}F_{\Gamma_t})-Ah$. We saw in \cite{lengeler} that in $S_\gamma$ we have
\begin{equation}\label{eqn:tildee}
 \tilde e-e=r(h/\gamma,hk,\nabla h)\circ\tau,
\end{equation}
where $r$ is an analytic function of its arguments. Thus, from
\[\phi_t^*(\grad_{L_2}F_{\Gamma_t})=\kappa\big(\Delta_{\tilde g} H_{\tilde e} + H_{\tilde e}(H_{\tilde e}^2/2-2K_{\tilde g})+C_0(2K_{\tilde g}-H_{\tilde e}C_0/2)\big)\]
and the results in \Appref{app:cov} we infer that $Q(h)$ is an analytic function of
zero to third order derivatives of $k$, and zero to fourth order derivatives of $h$. In \Secref{iteration} we will need the term $Q(h)$ to lie in $W^{1-1/p}_p(\Gamma)$. Since this term contains up to third order derivatives of $k$, assuming $\Gamma$ to be of class $W^{6-1/p}_p$ would in fact be sufficient.

We conclude that we have to analyze the following linear parabolic system
\begin{equation}\label{eqn:tildesystemlinear}
\begin{aligned}
\mu_b\Delta u-\grad\pi&=f_1&&\mbox{ in }\Omega\setminus\Gamma,\\
\Div u&=f_2&&\mbox{ in }\Omega\setminus\Gamma,\\
\mu\big(\Delta_{g} v + \grad_{g}(w\,H) + K\, v -2\Div_{g}(w\,k)\big)&\\
-\grad_{g} q+ P_\Gamma [\![S]\!]\nu&=f_3^\top&&\mbox{ on } \Gamma,\\
2\mu\big(\langle\nabla^{g} v,k\rangle_{g}
-w\,(H^2-2K)\big)- q\,H+[\![S]\!]\nu\cdot\nu-Ah&=f_3^\perp&&\mbox{ on }\Gamma,\\
\Div_{g} v- w\, H&=f_4&&\mbox{ on }\Gamma,\\
u-v-w\,\nu&=f_5&&\mbox{ on }\Gamma,\\
\partial_t h-w&=f_6&&\mbox{ on }\Gamma
\end{aligned}
\end{equation}
for suitable data $f_1,\ldots,f_6$ with the additional requirements $u=0$ on $\pa\Omega$ and $h(0)=h_0$ for some suitable initial value $h_0$. Here, we dropped the tilde symbols and, as before, $S=2\mu_b  Du - \pi I$.

\section{Linear Analysis}
\seclabel{linear}
In this section we study the linearization \eqref{eqn:tildesystemlinear} with fully inhomogeneous data and establish its unique solvability in the sense of maximal regularity in an $L_p$-setting. To begin with, let us specify suitable function spaces for the solutions and for the data. From \eqref{eqn:tildesystemlinear}$_{2,6}$ we obtain
\[\int_{\Gamma^i} (w+f_5\cdot\nu)\,dA=\int_{\Omega^i}f_2\,dx\]
for $i=1,\ldots,m$; combining this identity with \eqref{eqn:tildesystemlinear}$_5$ we obtain
\begin{equation}\label{comp1}
 \int_{\Gamma^i}f_4/H\,dA=-\int_{\Omega^i}f_2\,dx + \int_{\Gamma^i} f_5\cdot\nu\,dA \quad\text{for each $\Gamma^i$ that is a CMC surface}.
\end{equation}
Recall that the only closed, connected CMC (= constant mean curvature) surfaces embedded in $\R^3$ are the round spheres. Furthermore, of course, we have
\begin{equation}\label{comp2}
 \int_{\Omega}f_2\,dx=0.
\end{equation}
For $p\in (1,\infty)\setminus\{4\}$ we define
\begin{equation*}
\begin{aligned}
\G_p(T):=\big\{(f_1,\ldots,f_6,h_0)\,|\,f_1\in L_p(I,L_p(\Omega\setminus\Gamma,\R^3)),\ f_2\in L_p(I,H^1_p(\Omega\setminus\Gamma)),\\ f_3^\top\in L_p(I,L_p(\Gamma,T\Gamma)),\ f_3^\perp\in L_p(I,W^{1-1/p}_p(\Gamma)),\ f_4\in L_p(I,H^{1}_p(\Gamma)),\\
f_5\in L_p(I,W^{2-1/p}_p(\Gamma,\R^3)),\ f_6\in L_p(I,W^{2-1/p}_p(\Gamma)),\ h_0\in W^{5-4/p}_p(\Gamma)\big\},
\end{aligned}
\end{equation*}
where $I = (0,\,T)$, the space of data
\begin{equation*}
\begin{aligned}
\F_p(T):=\big\{(f_1,\ldots,f_6,h_0)\in \G_p(T)\,|\,\mbox{\eqref{comp1} and \eqref{comp2} hold for almost all $t$}\big\},
\end{aligned}
\end{equation*}
and the space of solutions
\begin{equation*}
\begin{aligned}
\E_p(T):=\big\{(u,v,w,\pi,q,h)\,|\,u\in L_p(I,H^2_p(\Omega\setminus\Gamma,\R^3)\cap {}_0 H^1_p(\Omega,\R^3)),\\ v\in L_p(I,H^2_p(\Gamma,T\Gamma)),\ w\in L_p(I,W^{2-1/p}_p(\Gamma)),\\ \pi\in L_p(I,H^1_p(\Omega\setminus\Gamma)),\ q\in L_p(I,H^1_p(\Gamma)), \\ h\in L_p(I,W^{5-1/p}_p(\Gamma))\cap H^1_p(I,W^{2-1/p}_p(\Gamma)),\\ \text{such that \eqref{comp1} and \eqref{comp2} with $f_2=\pi$, $f_4=q$,}\\ \text{and $f_5=0$ hold for almost all $t$}\big\};
\end{aligned}
\end{equation*}
each space is endowed with the canonical norm. 

\begin{theorem}\label{thm:linear}
For $p\in [2,\infty)\setminus\{4\}$ and $(f_1,\ldots,f_6,h_0)\in \F_p(T)$ there exists a unique solution $(u,v,w,\pi,q,h)\in \E_p(T)$ of \eqref{eqn:tildesystemlinear}. If the functions $f_1,\ldots,f_6$, and $h_0$ are smooth in space and time, then so is the solution $(u,v,w,\pi,q,h)$.\footnote{Here and in the following smoothness means $C^\infty$ up to possible jumps across $\Gamma$ for $f_1$, $f_2$, $\pi$, and first order derivatives of $u$.}
\end{theorem}

\proof[Proof: Existence for Smooth Data and Uniqueness] This follows by combining the elliptic theory proved in \cite{lengeler} with standard arguments from parabolic $L_2$-theory. 
We will successively eliminate the data $(f_1,\ldots,f_6)$ and hence write the velocity fields in the form $u=u_0+u_1+u_2$, $v=v_0+v_1+v_2$ (with strictly tangential $v_i$), and $w=w_0+w_1+w_2$. First, we eliminate $f_5$ and $f_6$ by choosing a smooth function $u_0$ such that $[u_0]_{\pa\Omega}=0$ and $[u_0]_\Gamma=f_5-f_6\,\nu$ and by defining $v_0:=0$ and $w_0:=-f_6$. Next, we eliminate $f_2$ and $f_4$ by solving the stationary system
\begin{equation*}
\begin{aligned}
\Div u_1&=f_2-\Div u_0&&\mbox{ in }\Omega\setminus\Gamma,\\
\DIV u_1&=f_4-\DIV u_0&&\mbox{ on }\Gamma
\end{aligned}
\end{equation*}
at fixed, but arbitrary $t\in\bar I$ for a smooth function $u_1$, see Theorem 3.6 in \cite{lengeler} with $f_1=0$ and $f_3=0$, and by choosing $v_1$, $w_1$ such that $u_1-v_1-w_1\nu=0$. Finally, we solve \eqref{eqn:tildesystemlinear} for $(u_2,v_2,w_2,\pi,q,h)$ with vanishing $f_2$, $f_4$, $f_5$, and $f_6$, with $f_1$ and $f_3$ replaced by $\tilde f_1:=f_1-2\mu_b\Div Du_1$ and $\tilde f_3:=f_3-2\mu\DIV\D u_1-[\![2\mu_b\,Du_1]\!]\nu$, respectively, and with $u_2-v_2-w_2\nu=0$. To this end, we note that this system can be written in the form
\begin{equation}\label{linnac}
\begin{aligned}
\Div S&=\tilde f_1 &&\mbox{ in }\Omega\setminus \Gamma,\\
\Div u_2&=0&&\mbox{ in }\Omega\setminus \Gamma,\\
\DIV {^fT}+[\![S]\!]\nu+Ah\,\nu&=\tilde f_3&&\mbox{ on }\Gamma,\\
\DIV u_2&=0&&\mbox{ on }\Gamma,\\
\pa_th-w_2&=0&&\mbox{ on }\Gamma,\\
\end{aligned}
\end{equation}
where the stress tensors $S$ and $T$ are taken with respect to $u_2$, $\pi$, and $q$. Let us multiply \eqref{linnac}$_3$ by a smooth test function $\varphi: \bar I\times \Omega\rightarrow \R^3$ fulfilling the divergence constraints \eqref{linnac}$_{2,4}$ and vanishing on $\pa\Omega$. Analogously to the computations in Section 3.1 of \cite{lengeler}, integration by parts then leads to the following weak formulation of our system
\begin{equation}\label{weak}
\begin{aligned}
B(u_2,\varphi)+A(h,\varphi) &= F(\varphi),\\
\pa_t h - w_2&=0,
\end{aligned}
\end{equation}
which is to hold for almost all $t\in I$. Here, we used the definitions
\begin{equation*}
\begin{aligned}
B(u,\varphi)&:=2\mu_b\int_{\Omega\setminus \Gamma} \langle Du,D\phi\rangle_e\,dx+2\mu\int_{\Gamma} \langle \D u,\D \phi\rangle_g\,dA,\\
A(h,\varphi)&:=\kappa\int_\Gamma \big(\Delta_g h\,\Delta_g \phi^\perp + a^{\alpha\beta}h_{,\alpha}\phi^\perp_{,\beta} + b\, h\,\phi^\perp\big)\,dA,\\
F(\varphi)&:=\int_\Gamma \langle \tilde f_3,\phi\rangle_e\,dA+\int_{\Omega\setminus \Gamma} \langle \tilde f_1,\phi\rangle_e\,dx
\end{aligned}
\end{equation*}
with $\phi^\perp:=\phi\cdot\nu$.
Choosing $\phi=u_2$ and making use of the coercivity of the bilinear form $B$, see Lemma 3.1 in \cite{lengeler}, and the $L^2$-theory of the Laplacian on $\Gamma$, it is not hard to see that we can estimate $u_2$ in $L_2(I,H^1_0(\Omega))$, $v_2$ in $L_2(I,H^1(\Gamma,T\Gamma))$, and $h$ in $L_\infty(I,H^2(\Gamma))\cap H^1(I,H^{1/2}(\Gamma))$ in terms of $\tilde f_3$ in $L_2(I,L_2(\Gamma))$, $\tilde f_1$ in $L_2(I,L_2(\Omega))$, and $h_0$ in $H^2(\Gamma)$. Thus, by Galerkin's method, see \cite{Evans:Partial-Differential-Equations}, we can actually construct such a \emph{weak solution} $(u_2,h)$. 
Next we reconstruct the pressure functions. Since, for fixed $t\in I$, the functional $\phi\mapsto B(u_2,\phi)+ A(h,\phi)-F(\phi)$ annihilates the space
\begin{equation*}
\begin{aligned}
X:=\Big\{\,u\in H^1_0(\Omega,\R^3)\,:\,\Div u=0 \text{ in }\Omega\setminus \Gamma,\ \DIV u=0 \text{ on } \Gamma,\ P_\Gamma u \in H^1(\Gamma;T\Gamma)\,\Big\},
\end{aligned}
\end{equation*}
by Corollary 3.3 in \cite{lengeler} there exist functions $(\pi,q)\in L_2(I,Z)$ with 
\[Z:=\Big\{\,(f_2,f_4)\in L_2(\Omega)\times L_2(\Gamma)\,:\,\text{\eqref{comp1}, \eqref{comp2} with $f_5=0$ hold}\,\Big\}\]
such that
\begin{equation*}
  B(u_2,\phi)+ A(h,\phi)-F(\phi)=-\int_{\Omega\setminus\Gamma}\pi\,\Div\phi\,dx - \int_{\Gamma}q\,\DIV\phi\,dA
 \end{equation*}
for all $\phi\in H_0^1(\Omega)$ with $P_\Gamma[\phi]_\Gamma\in H^1(\Gamma;T\Gamma)$ and almost all $t\in I$.
As announced above, the full (weak) solution of the system is then given by $(u,v,w,\pi,q)$, where $u:=u_0+u_1+u_2$, $v:=v_0+v_1+v_2$, and $w:=w_0+w_1+w_2$. Furthermore, we can estimate these functions in terms of the data analogously to the estimates (32) and (37) in \cite{lengeler}; this proves uniqueness in $\E_p(T)$ for $p\in [2,\infty)$ since then $\E_p(T)$ embeds into the above energy spaces.
%
It remains to prove smoothness of $(u_2,\pi,q)$. To this end we take the $k$-th derivative of \eqref{weak}$_1$ in time for arbitrary $k\in\N$ and choose $\phi=\pa_t^{k-1} u$; of course, strictly speaking this must be done on the Galerkin level. The resulting energy estimates show that $h\in H^{k-1}(I,H^2(\Omega))$; since $k$ is arbitrary, we have $h\in C^\infty(\bar I,H^2(\Gamma))$. Now, from Theorem 3.7 in \cite{lengeler} (with $f_5=\pa_t h$) we obtain $Ah\in C^\infty(\bar I,H^1(\Gamma))$. Thus, $L_2$-theory for the biharmonic operator on $\Gamma$ shows that $h\in C^\infty(\bar I,H^5(\Gamma))$; iterating this procedure we obtain that $h$ is smooth in space and time. Using Theorem 3.7 in \cite{lengeler} once more we see that the same is true for $u_2$, $\pi$, and $q$.
\qed
\medskip

The proof of existence in the case of non-smooth data in $\F_p(T)$, $p>2$, is of course much more involved. The main step is carried out in Subsection~\ref{linear:halfspace} where a maximal regularity result is shown for the principal linearization of our system in the prototype geometry of a double half-space.
In Subsection~\ref{sec:linear:domain} this result is then transfered to a bounded domain using the basic procedure of localization and transformation once more.

\subsection{Double half-space}
\label{linear:halfspace}
In this subsection we study the principal linearization of our system
\begin{equation}
	\eqnlabel{linear:halfspace}
	\begin{array}{rclll}
		                                                                 \eta\, u - \mu_b \Delta u + \grad \pi & = & f_1    & \quad \mbox{in} & \R^n\setminus\Sigma, \\[0.5em]
		                                                                                    \mbox{div}\,u & = & f_2    & \quad \mbox{in} & \R^n\setminus\Sigma, \\[0.5em]
		 - \mu \Delta v + \grad q - 2\mu_b P_\Sigma [\![Du]\!]\nu & = & f_3^\top & \quad \mbox{on} & \Sigma,      \\[0.5em]
		               - [\![\pi]\!] -\kappa \Delta^2 h & = & f_3^\perp  & \quad \mbox{on} & \Sigma,      \\[0.5em]
		                                                                    \Div v & = & f_4    & \quad \mbox{on} &\Sigma,      \\[0.5em]
		u - v-w\,\nu & = & f_5      & \quad \mbox{on} &\Sigma          \\[0.5em]                                                                                                                                                    (\partial_t+\eta) h - w & = & f_6      & \quad \mbox{on} &\Sigma	
	\end{array}
\end{equation}
with $h(0)=h_0$ posed in the unbounded time interval $\R_+ := (0,\infty)$ and in the prototype geometry $\R^n\setminus\Sigma$, where $n\ge 2$ and $\Sigma := \R^{n - 1}\times\{0\}$. We employ the splitting $(x,\,y) \in \R^{n - 1}\times \R$ for the spatial variables and $\nu:=e_y$. For technical reasons related to the unbounded spatial and temporal domain of the system we introduced an artificial shift $\eta > 0$. 
Note that due to $[\![u]\!] = 0$ we have $2 \mu_b [\![Du]\!] \nu \cdot \nu = 2 \mu_b [\![\partial_y u_n]\!] = 2 \mu_b [\![f_2]\!]$, which may be hidden in $f^\perp_3$.
For $p\in (1,\infty)$ we define the space of data
\begin{equation*}
\begin{aligned}
\F_p^\Sigma:=\big\{(f_1,\ldots,f_6,h_0)\,|\,f_1\in L_p(\R_+,L_p(\R^n\setminus\Sigma,\R^n)), f_2\in L_p(\R_+,H^1_p(\R^n\setminus\Sigma)),\\ f_3^\top\in L_p(\R_+,L_p(\Sigma,T\Sigma)),\ f_3^\perp\in L_p(\R_+,\dot W^{1-1/p}_p(\Sigma)),\ f_4\in L_p(\R_+,H^{1}_p(\Sigma)),\\
 f_5\in L_p(\R_+,W^{2-1/p}_p(\Sigma,\R^n)),\ f_6\in L_p(\R_+,W^{2-1/p}_p(\Sigma)),\ h_0\in W^{5-4/p}_p(\Sigma)\big\}
\end{aligned}
\end{equation*}
and the space of solutions
\begin{equation*}
\begin{aligned}
\E_p^\Sigma:=\big\{(u,v,w,\pi,q,h)\,|\,u\in L_p(\R_+,H^2_p(\R^n\setminus\Sigma,\R^n)\cap H^1_p(\R^n,\R^n)),\\ v\in L_p(\R_+,H^2_p(\Sigma,T\Sigma)),\ w\in L_p(\R_+,W^{2-1/p}_p(\Sigma)),\\ \pi\in L_p(\R_+,\dot{H}^1_p(\R^n\setminus\Sigma)/\R),\ q\in L_p(\R_+,\dot{H}^1_p(\Sigma)/\R),\\ h\in L_p(\R_+,W^{5-1/p}_p(\Sigma))\cap H^1_p(\R_+,W^{2-1/p}_p(\Sigma))\big\};
\end{aligned}
\end{equation*}
each space is endowed with the canonical norm.

\begin{theorem}
	\label{thm:linear:halfspace}
	For $\eta > 0$, $p\in(1,\infty)$, and $(f_1,\ldots,f_6,h_0)\in\F_p^\Sigma$ there exists a unique solution $(u,v,w,\pi,q,h)\in \E_p^\Sigma$ of \eqnref{linear:halfspace}.
\end{theorem}
\proof 
The proof will be carried out in two steps, where we split up the system into a stationary problem with inhomogeneous right-hand-sides and an evolution equation with homogeneous right hand sides. 

\subsection*{Step 1}
In this step we will eliminate all data except for $f_6$. To begin with, we eliminate $h_0$ by constructing an extension \[\bar h\in H^1_p(\R_+,\,W^{2 - 1/p}_p(\Sigma)) \cap L_p(\R_+,\,W^{5 - 1/p}_p(\Sigma));\] see the remark after Theorem \ref{thm}. In order to deal with the remaining data, we study the stationary system
\begin{equation*}
	\begin{array}{rclll}
		                                                                 \eta\, u - \mu_b \Delta u + \grad \pi & = & f_1    & \quad \mbox{in} & {\bR}^n\setminus\Sigma, \\[0.5em]
		                                                                                    \mbox{div}\,u & = & f_2    & \quad \mbox{in} & {\bR}^n\setminus\Sigma, \\[0.5em]
		                                                                           - \mu \Delta v + \grad q - 2\mu_b P_\Sigma [\![Du]\!]\nu & = & f_3^\top & \quad \mbox{on} & \Sigma,      \\[0.5em]
		                                                  - [\![\pi]\!] & = & f_3^\perp  & \quad \mbox{on} & \Sigma,      \\[0.5em]
 		                                                                    \Div v & = & f_4    & \quad \mbox{on} & \Sigma,\\[0.5em]
 		                                                                                                              u-v-w\,\nu& = & f_5    & \quad \mbox{on} & \Sigma.
	\end{array}
\end{equation*}
Concerning this system we show that there exists a unique solution 
\begin{equation*}
	\begin{array}{c}
		u \in H^2_p({\bR}^n\setminus\Sigma,\bR^n)\cap H^1_p(\R^n), \quad v\in H^2_p(\Sigma,T\Sigma), \quad w\in W^{2-1/p}_p(\Sigma), \\[0.5em]
		\pi \in \dot{H}^1_p({\bR}^n\setminus\Sigma)/\R, , \quad  q \in \dot{H}^1_p(\Sigma)/\R,
	\end{array}
\end{equation*}
provided that the data satisfy
\begin{equation*}
	\begin{array}{c}
		f_1\in L_p(\R^n\setminus\Sigma,\R^n),\quad f_2\in H^1_p(\R^n\setminus\Sigma),\quad f_3^\top\in L_p(\Sigma,T\Sigma),\\ f_3^\perp\in \dot W^{1-1/p}_p(\Sigma),\quad f_4\in H^{1}_p(\Sigma),\quad f_5\in W^{2-1/p}_p(\Sigma,\R^n).
	\end{array}
\end{equation*}
To begin with, we eliminate $f_5$ by choosing a function $\bar{u}\in H^2_p(\R^n\setminus\Sigma,\R^n)\cap H^1_p(\R^n)$ such that $[\bar{u}]_\Sigma=f_5$; for the surjectivity of the trace operator see, for instance, \cite{Triebel:Function-Spaces-1}.
With $f_5$ vanishing, we may employ the splitting $u=(v,w)\in \R^{n - 1}\times \R$ in the whole space $\R^n$. Next, we eliminate $f_1$, $f_2$, and $f_4$ by making use of the results of Appendix~\ref{sec:appendix:stokes}. To this end, we first solve the whole space problem
\begin{equation*}
	\begin{array}{rcll}
		\eta\, \bar{v} - \mu \Delta \bar{v} + \grad \bar{q} & = & 0   & \quad \mbox{on} \ \Sigma, \\[0.5em]
		                                      \Div\,\bar{v} & = & f_4 & \quad \mbox{on} \ \Sigma
	\end{array}
\end{equation*}
to obtain $\bar{v} \in H^2_p(\Sigma,\,T\Sigma)$ and $\bar{q} \in \dot{H}^1_p(\Sigma)$, and then we solve the two decoupled halfspace problems
\begin{equation*}
	\begin{array}{rcll}
		\eta\, \bar{u} - \mu_b \Delta \bar{u} + \grad \bar{\pi} & = & f_1     & \quad \mbox{in} \ \bR^n_\pm, \\[0.5em]
		                                   \mbox{div}\,\bar{u} & = & f_2     & \quad \mbox{in} \ \bR^n_\pm, \\[0.5em]
		                                    \bar{u} & = & \bar{v} & \quad \mbox{on} \ \Sigma
	\end{array}
\end{equation*}
to obtain $\bar u \in H^2_p(\R^n\setminus\Sigma,\R^n)\cap H^1_p(\R^n)$ and $\bar \pi \in \dot{H}^1_p(\R^n\setminus\Sigma)$.
Finally, in order to solve the reduced problem, we employ a Fourier transformation in the tangential variables to obtain the system
\begin{equation}
	\eqnlabel{linear:tpfourier}
	\begin{array}{rclll}
		                         \eta\, \hat{v} + \mu_b |\xi|^2 \hat{v} - \mu_b \partial^2_y \hat{v} + i \xi \hat{\pi} & = & 0,            & \quad \xi \in \bR^n, \ y \neq 0, \\[0.5em]
		                    \eta\, \hat{w} + \mu_b |\xi|^2 \hat{w} - \mu_b \partial^2_y \hat{v} + \partial_y \hat{\pi} & = & 0,            & \quad \xi \in \bR^n, \ y \neq 0, \\[0.5em]
		                                                                  i \xi^{\sfT}\,\hat{v} + \partial_y \hat{w} & = & 0,            & \quad \xi \in \bR^n, \ y \neq 0, \\[0.5em]
		                                                              [\![\hat{v}]\!] = 0, \quad [\![\hat{w}]\!] & = & 0             & \quad \xi \in \bR^n,\ y=0             \\[0.5em]
		\mu |\xi|^2 \hat{v} + i \xi \hat{q} - \mu_b [\![\partial_y \hat{v}]\!] - \mu_b i \xi [\![\hat{w}]\!] & = & \hat{g}_\tau, & \quad \xi \in \bR^n,\ y=0             \\[0.5em]
		                                                [\![\hat{\pi}]\!] & = & \hat{g}_\nu,  & \quad \xi \in \bR^n,\ y=0             \\[0.5em]
		                                                                                     i \xi^{\sfT}\hat{v} & = & 0,            & \quad \xi \in \bR^n,\ y=0
	\end{array}
\end{equation}
Here, we simplified the notation by setting $g = (g_\tau,\,g_\nu) := (f^\top_3,\,f^\bot_3)$.
The generic solution of the ODE system is easily seen to be given as
\begin{equation}\label{schwuff}
	\left[ \begin{array}{c} \hat{v}^\pm(\xi,\,y) \\[0.5em] \hat{w}^\pm(\xi,\,y) \\[0.5em] \hat{\pi}^\pm(\xi,\,y) \end{array} \right]
		= \left[ \begin{array}{rr} \varpi & \quad - i \zeta \\[0.5em] \pm i \zeta^{\sfT} & \quad \pm |\zeta| \\[0.5em] 0 & \quad \eta \sqrt{\mu_b} \end{array} \right]
		  \left[ \begin{array}{l} \hat{z}^\pm_v(\xi) e^{\mp \frac{\varpi}{\sqrt{\mu_b}} y} \\[0.5em] \hat{z}^\pm_w(\xi) e^{\mp |\xi| y} \end{array} \right], \quad \xi \in \bR^{n - 1}, \ y \gtrless 0
\end{equation}
with $\zeta := \sqrt{\mu_b}\,\xi$, $\varpi := \sqrt{\eta + |\zeta|^2}$, and four functions $\hat z^\pm_v: \bR^{n - 1} \rightarrow \bR^{n - 1}$, and $\hat z^\pm_w: \bR^{n - 1} \rightarrow \bR$,
which have to be determined together with $q: \bR^{n - 1} \rightarrow \bR$ based on the \emph{transmission conditions} \eqnref{linear:tpfourier}$_{4,5,6}$ and the incompressibility constraint on the membrane. With this representation of the solution the transmission conditions become
\begin{equation*}
	\begin{array}{c}
		\varpi (\hat{z}^+_v - \hat{z}^-_v) - i \zeta (\hat{z}^+_w - \hat{z}^-_w) = 0, \qquad i \zeta^{\sfT} (\hat{z}^+_v + \hat{z}^-_v) + |\zeta| (\hat{z}^+_w + \hat{z}^-_w) = 0 \\[0.5em]
		\frac{\mu}{\mu_b} |\zeta|^2 (\varpi \hat{z}^+_v - i \zeta \hat{z}^+_w) + \frac{1}{\sqrt{\mu_b}} i \zeta \hat{q} + \sqrt{\mu_b} \varpi^2 (\hat{z}^+_v + \hat{z}^-_v) - \sqrt{\mu_b} i \zeta |\zeta| (\hat{z}^+_w + \hat{z}^-_w) = \hat{g}_\tau \\[0.5em]
		2 \sqrt{\mu_b} \varpi i \zeta^{\sfT} (\hat{z}^+_v - \hat{z}^-_v) + 2 \sqrt{\mu_b} |\zeta|^2 (\hat{z}^+_w - \hat{z}^-_w) + \eta \sqrt{\mu_b} (\hat{z}^+_w - \hat{z}^-_w) = \hat{g}_\nu,
	\end{array}
\end{equation*}
and the incompressibility constraint on the membrane reads
\begin{equation*}
	i \zeta^{\sfT} (\varpi \hat{z}^+_v - i \zeta \hat{z}^+_w) = i \zeta^{\sfT} (\varpi \hat{z}^-_v - i \zeta \hat{z}^-_w) = 0.
\end{equation*}
Applying $i \zeta^{\sfT}\,\cdot$ to the tangential transmission condition and using the continuity of $w$ across the membrane and the incompressibility condition we obtain
\begin{equation*}
	- \frac{1}{\sqrt{\mu_b}} |\zeta|^2 \hat{q} - \eta \sqrt{\mu_b} |\zeta| (\hat{z}^+_w + \hat{z}^-_w) = i \zeta^{\sfT} \hat{g}_\tau,
\end{equation*}
which leads to
\begin{equation}\label{qglg}
	\frac{1}{\sqrt{\mu_b}} i \zeta \hat{q} = - \eta \sqrt{\mu_b} \frac{i \zeta}{|\zeta|} (\hat{z}^+_w + \hat{z}^-_w) - \frac{i \zeta \otimes i \zeta}{|\zeta|^2} \hat{g}_\tau.
\end{equation}
For the tangential transmission condition we then obtain
\begin{equation*}
	\begin{array}{r}
		{\displaystyle \frac{\mu}{\mu_b} |\zeta|^2 (\varpi \hat{z}^+_v - i \zeta \hat{z}^+_w)}
			+ {\displaystyle \sqrt{\mu_b} \varpi^2 (\hat{z}^+_v + \hat{z}^-_v)}
			- {\displaystyle \eta \sqrt{\mu_b} \frac{i \zeta}{|\zeta|} (\hat{z}^+_w + \hat{z}^-_w)} \qquad \qquad \qquad \\[1.5em]
		= {\displaystyle \left( 1 + \frac{i \zeta \otimes i \zeta}{|\zeta|^2} \right) \hat{g}_\tau.}
	\end{array}
\end{equation*}
Since the continuity of $v$ across the membrane implies
\begin{equation*}
	\frac{\mu}{\mu_b} |\zeta|^2 (\varpi \hat{z}^+_v - i \zeta \hat{z}^+_w)
		= \frac{1}{2} \frac{\mu}{\mu_b} |\zeta|^2 \Big( \varpi (\hat{z}^+_v + \hat{z}^-_v) - i \zeta (\hat{z}^+_w + \hat{z}^-_w) \Big),
\end{equation*}
the tangential transmission condition may be rewritten as
\begin{equation*}
	\begin{array}{r}
		{\displaystyle \left( \sqrt{\mu_b} \varpi + \frac{1}{2} \frac{\mu}{\mu_b} |\zeta|^2 \right) \varpi (\hat{z}^+_v + \hat{z}^-_v)}
			- {\displaystyle \left( \frac{\eta \sqrt{\mu_b}}{|\zeta|} - \frac{1}{2} \frac{\mu}{\mu_b} |\zeta|^2 \right) i \zeta (\hat{z}^+_w + \hat{z}^-_w)} \qquad \qquad \\[1.5em]
		= {\displaystyle \left( 1 + \frac{i \zeta \otimes i \zeta}{|\zeta|^2} \right) \hat{g}_\tau.}
	\end{array}
\end{equation*}
Furthermore, due to the continuity of $v$ across the membrane, the normal transmission condition simplifies to
\begin{equation*}
	\eta \sqrt{\mu_b} (\hat{z}^+_w - \hat{z}^-_w) = \hat{g}_\nu.
\end{equation*}
On the other hand, the continuity of $w$ across the membrane together with the incompressibility constraint on the membrane, which may equivalently be written in the form
\begin{equation*}
	\varpi i \zeta^{\sfT} (\hat{z}^+_v + \hat{z}^-_v) = - |\zeta|^2 (\hat{z}^+_w + \hat{z}^-_w),
\end{equation*}
imply
\begin{equation*}
	\varpi |\zeta| (\hat{z}^+_w + \hat{z}^-_w) = - \varpi i \zeta^{\sfT} (\hat{z}^+_v + \hat{z}^-_v) = |\zeta|^2 (\hat{z}^+_w + \hat{z}^-_w).
\end{equation*}
This yields $(\varpi - |\zeta|) |\zeta| (\hat{z}^+_w + \hat{z}^-_w) = 0$, that is, $\hat{z}^+_w + \hat{z}^-_w = 0$ and thus  $i \zeta^{\sfT}(\hat{z}^+_v + \hat{z}^-_v) = 0$. Hence, we obtain
\begin{equation}
	\eqnlabel{tpzw}
	\hat{z}^\pm_w = \pm \frac{1}{2}\,\frac{1}{\eta \sqrt{\mu_b}}\,\hat{g}_\nu.
\end{equation}
Combining these identities with the tangential transmission condition and the continuity of $v$ across the membrane we infer
\begin{equation*}
	\begin{array}{rcl}
		{\displaystyle \left( \sqrt{\mu_b} \varpi + \frac{1}{2} \frac{\mu}{\mu_b} |\zeta|^2 \right) \varpi (\hat{z}^+_v + \hat{z}^-_v)}
			& = & {\displaystyle \left( 1 + \frac{i \zeta \otimes i \zeta}{|\zeta|^2} \right) \hat{g}_\tau}, \\[1.5em]
		\varpi (\hat{z}^+_v - \hat{z}^-_v)
			& = & {\displaystyle \frac{i \zeta}{\eta \sqrt{\mu_b}}\,\hat{g}_\nu}.
	\end{array}
\end{equation*}
Adding and subtracting these two equations yields
\begin{equation}
	\eqnlabel{tpzv}
	\varpi \hat{z}^\pm_v = \frac{1}{2} \left( \sqrt{\mu_b} \varpi + \frac{1}{2} \frac{\mu}{\mu_b} |\zeta|^2 \right)^{-1} \left( 1 + \frac{i \zeta \otimes i \zeta}{|\zeta|^2} \right) \hat{g}_\tau \pm \frac{1}{2}\,\frac{i \zeta}{\eta \sqrt{\mu_b}} \hat{g}_\nu.
\end{equation}
Combining \eqnref{tpzw} and \eqnref{tpzv} we find
\begin{equation*}
	[\hat{v}]_\Sigma = \varpi \hat{z}^\pm_v - i \zeta \hat{z}^\pm_w = \frac{1}{2} \left( \sqrt{\mu_b} \varpi + \frac{1}{2} \frac{\mu}{\mu_b} |\zeta|^2 \right)^{-1} \left( 1 + \frac{i \zeta \otimes i \zeta}{|\zeta|^2} \right) \hat{g}_\tau.
\end{equation*}
Now, the last symbol on the right-hand-side belongs to the Helmholtz projection
\begin{equation*}
	\cH_\Sigma: L_p(\Sigma,T\Sigma) \rightarrow L_{p, \sigma}(\Sigma,T\Sigma),
\end{equation*}
the projection associated to the direct topological decomposition
\begin{equation*}
	L_p(\Sigma,T\Sigma) = L_{p, \sigma}(\Sigma,T\Sigma) \oplus \nabla \dot{H}^1_p(\Sigma,T\Sigma),
\end{equation*}
where $L_{p,\sigma}(\Sigma,T\Sigma)\subset L_{p}(\Sigma,T\Sigma)$ denotes the subspace of solenoidal vector fields. Observe that $\cH_\Sigma$ is bounded as follows for instance from Mikhlin's multiplier theorem.
Based on this observation we may write
\begin{equation*}
	(\mu_b + |\zeta|^2) [\hat{v}]_\Sigma = \frac{1}{2} \frac{\mu_b + |\zeta|^2}{\sqrt{\mu_b} \varpi + \frac{1}{2} \frac{\mu}{\mu_b} |\zeta|^2} \widehat{\cH_\Sigma g_\tau}
\end{equation*}
and infer that $[v]_\Sigma \in H^2_p(\Sigma,T\Sigma)$ from $g_\tau \in L_p(\Sigma,T\Sigma)$, as follows again from Mikhlin's multiplier theorem and the characterization of Sobolev spaces via Bessel potentials;
see, for instance, the theorem of Section~2.5.6 in \cite{Triebel:Function-Spaces-1}.
Now, \eqref{qglg} simplifies to
\begin{equation*}
	- \frac{1}{\sqrt{\mu_b}} i \zeta \hat{q} =  \frac{i \zeta \otimes i \zeta}{|\zeta|^2} \hat{g}_\tau,
\end{equation*}
which yields $\grad q \in L_p(\Sigma,T\Sigma)$.
Finally, we have
\begin{equation*}
	[\hat{\pi}^\pm]_\Sigma = \eta \sqrt{\mu_b} \hat{z}^\pm_w = \pm \frac{1}{2} \hat{g}_\nu,
\end{equation*}
which yields $\bar g_\nu^\pm:=[\pi^\pm]_\Sigma \in \dot{W}^{1 - 1/p}_p(\Sigma)$, and since the solution constructed above also satisfies
the two decoupled Stokes systems
\begin{equation*}
	\begin{array}{rcll}
		\eta\, u - \mu_b \Delta u + \grad \pi & = & 0               & \quad \mbox{in} \ \bR^n_\pm, \\[0.5em]
		                       \mbox{div}\,u & = & 0               & \quad \mbox{in} \ \bR^n_\pm, \\[0.5em]
		                                   v & = & \bar{g}_\tau    & \quad \mbox{on} \ \Sigma,    \\[0.5em]
		                                 \pi & = & \bar{g}^\pm_\nu & \quad \mbox{on} \ \Sigma,
	\end{array}
\end{equation*}
with $\bar{g}_\tau:=[v]_\Sigma$, we obtain the desired regularity for $u$ and $\pi$; see Appendix~\ref{sec:appendix:stokes}. Note that the computations above imply
\begin{equation}
	\eqnlabel{linear:tptrace}
	[\hat{w}]_\Sigma = \frac{1}{2}\,\frac{|\zeta|}{\sqrt{\mu_b} \varpi (\varpi + |\zeta|)}\,\hat{g}_\nu;
\end{equation}
in particular, the trace of the normal component of the velocity field depends only on the right hand side of the normal transmission condition.

\subsection*{Step 2}
Next, we employ a Laplace transform in the time variable and a Fourier transform in the tangential space variables in order to compute the boundary symbol of the reduced problem, which will then be used to derive the exact mapping properties of the solution operator $f \mapsto h$; here, we simplified the notation by setting $f := f_6$. Thus, we consider the transformed system
\begin{equation}
	\eqnlabel{linear:fourier}
	\begin{array}{rcll}
		                         \eta\, \hat{v} + \mu_b |\xi|^2 \hat{v} - \mu_b \partial^2_y \hat{v} + i \xi \hat{\pi} & = & 0,       & \quad \lambda \in \Sigma_\theta,\ \xi \in \bR^n,\ y \neq 0, \\[0.5em]
		                    \eta\, \hat{w} + \mu_b |\xi|^2 \hat{w} - \mu_b \partial^2_y \hat{w} + \partial_y \hat{\pi} & = & 0,       & \quad \lambda \in \Sigma_\theta,\ \xi \in \bR^n,\ y \neq 0, \\[0.5em]
		                                                                     i \xi^\sfT \hat{v} + \partial_y \hat{w} & = & 0,       & \quad \lambda \in \Sigma_\theta,\ \xi \in \bR^n,\ y \neq 0, \\[0.5em]
		                                                              [\![\hat{v}]\!] = 0, \quad [\![\hat{w}]\!] & = & 0,       & \quad \lambda \in \Sigma_\theta,\ \xi \in \bR^n,\ y = 0, \\[0.5em]
		\mu |\xi|^2 \hat{v} + i \xi \hat{q} - \mu_b [\![\partial_y \hat{v}]\!] - \mu_b i \xi [\![\hat{w}]\!] & = & 0,       & \quad \lambda \in \Sigma_\theta,\ \xi \in \bR^n,\ y = 0,            \\[0.5em]
		                         \kappa |\xi|^4 \hat{h}  + [\![\hat{\pi}]\!] & = & 0,       & \quad \lambda \in \Sigma_\theta,\ \xi \in \bR^n,           \ y = 0, \\[0.5em]
		                                                                                     i \xi^{\sfT}[\hat{v}]_y & = & 0,       & \quad \lambda \in \Sigma_\theta,\ \xi \in \bR^n, \ y = 0,           \\[0.5em]
		                                                                          \lambda_\eta\, \hat{h} - [\hat{w}]_y & = & \hat{f}, & \quad \lambda \in \Sigma_\theta,\ \xi \in \bR^n,\ y = 0,
	\end{array}
\end{equation}
where we employ the abbreviation $\lambda_\eta := \lambda+\eta$ and denote by
\begin{equation*}
	\Sigma_\theta := \{\,z \in \bC\,:\,z \neq 0,\ |\mbox{arg}\,z| < \theta\,\}
\end{equation*}
a sector in the complex plane with opening angle $\frac{\pi}{2} < \theta < \pi$. 
To solve the transformed system \eqnref{linear:fourier} we reuse the computations made in the first step
and consider the first seven lines as an instance of problem \eqnref{linear:tpfourier} with right hand sides $\hat{g}_\tau = 0$ and $\hat{g}_\nu = -\frac{\kappa}{\mu^2_b} |\zeta|^4 \hat{h}$.
Then formula \eqnref{linear:tptrace} yields
\begin{equation*}
	[\hat{w}]_\Sigma = \frac{1}{2}\,\frac{|\zeta|}{\sqrt{\mu_b} \varpi (\varpi + |\zeta|)}\,\hat{g}_\nu
		= - \alpha\,\frac{|\zeta|}{\varpi (\varpi + |\zeta|)}\,|\zeta|^4 \hat{h}
\end{equation*}
with $\alpha := \frac{1}{2} \kappa / \mu^{5 / 2}_b > 0$, and we obtain
\begin{equation*}
	s(\lambda,\,|\xi|) \hat{h} := \left( \lambda_\eta + \alpha \frac{|\zeta|}{\varpi (\varpi + |\zeta|)}\,|\zeta|^4 \right) \hat{h} = \hat{f}.
\end{equation*}
Obviously, the boundary symbol $s$ has no zeros, if $\lambda \in \Sigma_\theta$ with $0 \leq \theta < \pi$.
Thus, the equation $s \hat{h} = \hat{f}$ may be uniquely solved for $\hat{h}$ and problem \eqnref{linear:halfspace} admits a unique solution - at least in the sense of tempered distributions. To prove the regularity assertions on $h$ a more precise analysis of the boundary symbol $s$ is necessary. To this end, we now consider the complex symbol
\begin{equation*}
	s(\lambda,\,z) = \lambda_\eta + m(z)n(z)\quad \mbox{with}
		\quad m(z) = \alpha \frac{\varpi(z)}{\varpi(z) + z},
		\ n(z) = \frac{z^5}{\varpi(z)^2},
\end{equation*}
where $\varpi(z) := \sqrt{\eta + z^2}$, $z \in \Sigma_\vartheta$ with $0 \leq \vartheta < \frac{\pi}{2}$, and $\lambda \in \Sigma_\theta$. Note that $|m(z)|$ is uniformly positive and bounded on $\bar{\Sigma}_\vartheta$; in particular, we have
\begin{equation*}
	|m(z)n(z)| \geq c(\vartheta)\,|n(z)|
\end{equation*}
for all $z \in \Sigma_\vartheta$ and some constant $c(\vartheta) > 0$. Moreover, note that $\lambda_\eta\in\Sigma_\theta$ for $\lambda\in\Sigma_\theta$ as well as $m(z)\in\Sigma_{2\vartheta}$, $n(z)\in\Sigma_{7\vartheta}$ for $z\in\Sigma_{\vartheta}$. Hence, we can easily prove by contradiction that, assuming $0 < 9 \vartheta < \pi - \theta$, we have
\begin{equation*}
	|\lambda_\eta + m(z)\,n(z)| \geq c(\theta,\,\vartheta)\,(|\lambda_\eta| + |m(z)n(z)|)
\end{equation*}
for all $z \in \Sigma_\vartheta$, $\lambda\in\Sigma_\theta$, and some constant $c(\theta,\vartheta) > 0$. These estimates imply that for the functions
\begin{equation*}
	(\lambda,\,z) \mapsto \lambda_\eta / s(\lambda,\,z) =: \phi(\lambda,\,z),
		\qquad (\lambda,\,z) \mapsto n(z) / s(\lambda,\,z) =: \psi(\lambda,\,z)
\end{equation*}
we have
\begin{equation}
	\label{holomorphic-symbols}
	\phi \in \cH^\infty(\Sigma_\theta\times\Sigma_\vartheta)\quad \mbox{ and }\quad \psi \in \cH^\infty(\Sigma_\theta\times\Sigma_\vartheta),
\end{equation}
provided that $\pi/2 < \theta < \pi$ and $0 < \vartheta < (\pi - \theta)/9$,
where we denote by $\cH^\infty$ the spaces of bounded holomorphic functions.
The desired regularity of $h$ may now be obtained as follows:
First observe that the operators
\begin{equation*}
	\begin{array}{rrcl}
		      \partial_t: & {}_0 H^1_p(\bR_+,\,W^{2 - 1/p}_p(\Sigma)) & \subseteq & X \longrightarrow X, \\[0.5em]
		(- \Delta)^{1/2}: &        L_p(\bR_+,\,W^{3 - 1/p}_p(\Sigma)) & \subseteq & X \longrightarrow X
	\end{array}
\end{equation*}
admit bounded $\cH^\infty$-calculi in the space $X:=L_p(\bR_+,\,W^{2 - 1/p}_p(\Sigma))$ with angles $\alpha^\infty_{\partial_t} = \frac{\pi}{2}$ and $\alpha^\infty_{(- \Delta)^{1/2}} = 0$, that is, these operators admit functional calculi
\begin{equation*}
	\altphi \mapsto \altphi(\partial_t): \cH^\infty(\Sigma_\theta) \rightarrow \cB(X), \quad
	\altphi \mapsto \altphi((- \Delta)^{1/2}): \cH^\infty(\Sigma_\vartheta) \rightarrow \cB(X)
\end{equation*}
provided that $\alpha^\infty_{\partial_t} < \theta < \pi$ and $\alpha^\infty_{(- \Delta)^{1/2}} < \vartheta < \pi$; see, for instance, Corollary 2.10 in \cite{Denk-Saal-Seiler:Newton-Polygon}. Moreover, the same corollary shows that we may employ Theorem~6.1 in \cite{Kalton-Weis:Operator-Sums} to obtain a joint $\cH^\infty$-calculus for these operators, that is, a functional calculus
\begin{equation*}
	\altphi \mapsto \altphi(\partial_t,\,(- \Delta)^{1/2}): \cH^\infty(\Sigma_\theta \times \Sigma_\vartheta) \rightarrow \cB(X)
\end{equation*}
provided that $\alpha^\infty_{\partial_t} < \theta < \pi$ and $\alpha^\infty_{(- \Delta)^{1/2}} < \vartheta < \pi$. It is shown, for instance, in \cite{Denk-Saal-Seiler:Newton-Polygon} that the operators $\altphi(\partial_t,\,(- \Delta)^{1/2})$ are Fourier-Laplace multipliers whose symbols are given by $\altphi(\lambda,\,|\xi|)$. Therefore, due to \eqref{holomorphic-symbols} and this joint $\cH^\infty$-calculus we infer that
\begin{equation*}
	\begin{array}{r}
		(\eta + \partial_t) h = \phi(\partial_t,\,(-\Delta)^{1/2})\,f \\[0.5em]
		(-\mu_b\Delta)^{5/2} (\eta - \mu_b\Delta)^{-1} h = \psi(\partial_t,\,(-\Delta)^{1/2})\,f
	\end{array}
	\quad \in L_p(\R_+,\,W^{2 - 1/p}_p(\Sigma)),
\end{equation*}
which implies that $h$ belongs to the asserted regularity class. This completes the proof of Theorem \ref{thm:linear:halfspace}. \qed
\medskip

\subsection{Bounded domain}
\seclabel{linear:domain}
Let us now finish the proof of Theorem \ref{thm:linear}. In view of the smoothness and uniqueness part of this theorem (which we already showed), by density it is sufficient to prove the estimate
\begin{equation*}
 \|(u,v,w,\pi,q,h)\|_{\E_p(T)} \le c\|(f_1,\ldots,f_6,h)\|_{\F_p(T)}
\end{equation*}
for smooth data. This estimate can be reduced to the assertion of Theorem \ref{thm:linear:halfspace} by the classical techniques of localization and transformation. We will only give a brief sketch of the procedure; see also the proof of Theorem 3.6 in \cite{lengeler}. To begin with, we note that in fact it is sufficient to prove the inequality
\begin{equation}\label{vollest2}
\begin{aligned}
\|(u,v,w,\pi,q,h)\|_{\E_p(T)} \le c\big(&\|(f_1,\ldots,f_6,h)\|_{\F_p(T)}
+\|\nabla u\|_{L_p(I\times\Omega)}\\& +\|\pi\|_{L_p(I\times\Omega)}+\|q\|_{L_p(I\times\Gamma)}+\|h\|_{L_p(I\times\Gamma)}\big).
\end{aligned}
\end{equation}
Indeed, combining this estimate with the uniqueness of solutions in $\E_p(T)$, $p\ge 2$, a standard contradiction argument shows that
\[\|\nabla u\|_{L_p(I\times\Omega)}+\|\pi\|_{L_p(I\times\Omega)}+\|q\|_{L_p(I\times\Gamma)}+\|h\|_{L_p(I\times\Gamma)}\le c\|(f_1,\ldots,f_6,h)\|_{\F_p(T)}.\]
The next step is to see that we can assume without restriction the solution to be localized in space. Indeed, if this is not the case we can multiply the solution by finitely many smooth cut-off functions; each of the products then solves the system \eqref{eqn:tildesystemlinear} where the right hand sides $f_1,\ldots,f_6$ now contain additional expressions involving lower order derivatives of the solution. Combining these finitely many estimates and using interpolation and absorption we arrive at \eqref{vollest2}. Now, if the spatial support of our the solution is strictly contained in $\bar\Omega\setminus\Gamma$ we can use standard results from $L_p$-theory of the Stokes system (see for instance \cite{Galdi:Navier-Stokes-1}) to prove \eqref{vollest2}. On the other hand, if the spatial support intersects $\Gamma$ we have to reduce the problem to Theorem \ref{thm:linear:halfspace}. In this case let us assume that the solution is supported in an open cube $Q_R$ of side length $R>0$ which is centered at some point $x_0\in\Gamma$. Rotating and translating the Cartesian coordinate system and choosing $R$ smaller if necessary, we may assume that $x_0=0$ and that $\Gamma\cap Q_R$ is the graph of a 
smooth function $a: Q_R^2:=Q_R\cap \Sigma\rightarrow (-R/2,R/2)$ such that $a(0)=0$ and $\nabla a(0)=0$. Consider the smooth diffeomorphism 
\[\Phi^{-1}: Q_R\rightarrow \tilde  Q_R:=\Phi^{-1}(Q_R), (x',x^3)\mapsto (x',x^3-a(x')).\]
This diffeomorphism induces the metric $\tilde e:=\Phi^*e$ on $\tilde Q_R$. We denote the restriction of $\tilde e$ to $Q_R^2$ by $\tilde g$. Note that $\Phi:(\tilde Q_R,\tilde e)\rightarrow (Q_R,e)$ and $\Phi|_{Q_R^2}:(Q_R^2,\tilde g)\rightarrow (\Gamma\cap Q_R,g)$ are isometries. Let us denote the pullbacks of the involved fields by $\tilde u:=\Phi^*u$, $\tilde\pi:=\Phi^*\pi$, $\tilde v:=\Phi^*v$, $\tilde w:=\Phi^*w$, $\tilde q:=\Phi^*q$, $\tilde h:=\Phi^*h$, $\tilde f_3^\top:=\Phi^*(P_\Gamma f_3)$, $\tilde f_3^\perp:=\Phi^*(f_3\cdot\nu)$, $\tilde f_i=\Phi^*f_i$ for $i=1,2,4,5,6$, and $\tilde h_0=\Phi^*h_0$. Now, proceeding as in \Secref{hanzawa}, that is, exploiting naturality of covariant differentiation under isometries and using the results from \Appref{app:cov} we see that  \eqref{eqn:tildesystemlinear} can be written in the form
\begin{equation*}
\begin{aligned}
-\eta\,\tilde u+\mu_b\Delta\tilde u-\grad\tilde\pi&=\hat f_1&&\mbox{ in }\R^3\setminus\Sigma,\\
\Div \tilde u&=\hat f_2&&\mbox{ in }\R^3\setminus\Sigma,\\
\mu\Delta \tilde v -\grad \tilde q  +2\mu_b[\![D\tilde u]\!]\nu&=\hat f_3^\top&&\mbox{ on } \Sigma,\\
-[\![\tilde\pi]\!]-\kappa\Delta^2\tilde h&=\hat f_3^\perp&&\mbox{ on }\Sigma,\\
\Div \tilde v&=\hat f_4&&\mbox{ on }\Sigma,\\
\tilde u-\tilde v-\tilde w\,\nu&=\hat f_5&&\mbox{ on }\Sigma,\\
(\partial_t+\eta) \tilde h-\tilde w&=\hat f_6&&\mbox{ on }\Sigma
\end{aligned}
\end{equation*}
with
\begin{equation*}
\begin{aligned}
\hat f_1&=\tilde f_1 + (\tilde e-e)*r(\tilde e)*(\mu_b\nabla^2 \tilde u,\grad\tilde\pi) - \eta\,\tilde u\\&\quad + \mu_b\, r(\tilde e)*\big((\nabla^2\tilde e,(\nabla\tilde e)^2)*\tilde u+\nabla\tilde e*\nabla \tilde u\big),\\
\hat f_2&=\tilde f_2+r(\tilde e)*\nabla\tilde e*\tilde u,\\
\hat f_3^\top
&=\tilde f_3^\top+(\tilde e-e)*r(\tilde e)*(\mu(\nabla^g)^2 \tilde v,\grad_g \tilde q) + \mu_b\,r(\tilde e)*\big([\nabla \tilde u]+\nabla\tilde e*[\tilde u]\big)\\
&\quad+ \mu\,r(\tilde e)*\big((\nabla^2\tilde e,(\nabla\tilde e)^2)*[\tilde u] + \nabla\tilde e*[\nabla\tilde u]\big),\\
\hat f_3^\perp&=\tilde f_3^\perp + \mu\,r(\tilde e)*\big(\nabla\tilde e*\nabla^g\tilde v + (\nabla\tilde e)^2*[\tilde u]\big)
+ r(\tilde e)*\nabla\tilde e\,\tilde q + (\tilde e-e)*r(\tilde e)*\nabla^4h\\ &\quad + \mbox{terms depending linearly on up to third order derivatives of }h,\\
\hat f_4&=\tilde f_4 + r(\tilde e)*\nabla\tilde e*[\tilde u],\\
\hat f_5&=\tilde f_5 + (\tilde e-e)*r(\tilde e)\,\tilde w,\\
\hat f_6&=\tilde f_6 + (\tilde e-e)*r(\tilde e)\,\pa_t\tilde h + \eta\,h.
\end{aligned}
\end{equation*}
Furthermore, it is not hard to see that
\[\tilde e(x',x^3)-e=r(\nabla a(x'))\]
with an analytic function $r$ such that $r(0)=0$. Theorem \ref{thm:linear:halfspace} and the open mapping theorem show that there exists a constant $c>0$ such that
\[\|(\tilde u,\tilde v,\tilde w,\tilde \pi,\tilde q,\tilde h)\|_{\E_p^\Sigma}\le c\|(\hat f_1,\ldots,\hat f_6,h_0)\|_{\F_p^\Sigma},\]
where the data $\hat f_1,\ldots,\hat f_6$ is extended to $\R_+$ by $0$ and  $(\tilde u,\tilde v,\tilde w,\tilde \pi,\tilde q,\tilde h)$ denotes the unique continuation of our solution to $\R_+$ which exists according to Theorem \ref{thm:linear:halfspace}. Making $\|\tilde e-e\|_{L_\infty(\tilde Q_R)}$ sufficiently small (by choosing $R$ small) for the highest order terms in $\hat f_1,\ldots,\hat f_6$ and using interpolation and Young's inequality for the lower order terms, by absorption we obtain
\begin{equation*}
\begin{aligned}
\|(\tilde u,\tilde v,\tilde w,\tilde \pi,\tilde q,\tilde h)\|_{\E_p^\Sigma} \le c\big(&\|(\tilde f_1,\ldots,\tilde f_6,\tilde h_0)\|_{\F_p(T)}
+\|\nabla \tilde u\|_{L_p(I\times\R^3)}\\& +\|\tilde \pi\|_{L_p(I\times\R^3)}+\|\tilde q\|_{L_p(I\times\Sigma)}+\|\tilde h\|_{L_p(I\times\Sigma)}\big).
\end{aligned}
\end{equation*}
Transforming this estimate back to $\Omega$ and $\Gamma$ and using once more interpolation and absorption to deal with the lower terms arising we arrive at \eqref{vollest2}. We omit the details. This proves Theorem \ref{thm:linear}.
\qed

\section{Contraction}
\seclabel{iteration}
In this section we finish the proof of the main result. For $1 < p < \infty $ with $p \neq 4$ we define
\[ \E_p^6(T):=L_p(I,\,W^{5 - 1/p}_p(\Gamma)) \cap H^1_p(I,\,W^{2 - 1/p}_p(\Gamma)).\]
Then the following embeddings are valid.

\begin{lemma}\label{lemma:emb}
 For $1 < p < \infty $ with $p \neq 4$ we have 
 \begin{itemize}
  \item[(i)] $\E_p^6(T)\hookrightarrow C(\bar I,W^{5-4/p}_p(\Gamma))\hookrightarrow C(\bar I,C^3(\Gamma))$,
  \item[(ii)] $\big\{h\in \E_p^6(T)\,|\,h(0)=0\big\}\hookrightarrow C(\bar I,W^{5-4/p}_p(\Gamma))\hookrightarrow C(\bar I,C^3(\Gamma))$, where the embedding constants are independent of $T$.
 \end{itemize}
\end{lemma}
\proof The embedding (i) follows from Theorem 4.10.2 in Chapter III of \cite{amann95} and the theorem in Section 7.4.4 of \cite{Triebel:Function-Spaces-2}; obviously, the embedding constant remains uniformly bounded as long as $T>0$ is bounded from below. The second embedding is a consequence of Remark 2 in Section 2.7.1 of \cite{Triebel:Function-Spaces-1} and a localization procedure; cf. \cite{lengeler}.
Now, (ii) follows from (i) by extending $h$ to the negative half-line by $0$.
\qed
\medskip

We denote by $L$ the linear parabolic operator defined by the left-hand-side of \eqref{eqn:tildesystem2}, and we consider $N:=(N_1,\ldots,N_6)$ as a nonlinear function of $(u,v,w,\pi,q,h)$. For $\delta>0$ let
\[C_\delta(T):=\Big\{\,(u,v,w,\pi,q,h)\in \bar B_\delta(0)\subset\E_p(T)\,:\,\|h\|_{L_\infty((0,T)\times\Gamma)}\le \gamma/2\,\Big\}.\]
Then, the function $\beta$ in the construction of $\Phi_h$, $h\in \cup_{\delta>0}C_\delta$, in the beginning of \Secref{hanzawa} can be chosen to be fixed, and, in particular, the generic analytic funtions $r$ in the nonlinearities do not dependend on $h$. For $1 < p < \infty$ let
\[{_0\E_p}(T):=\big\{\,(u,v,w,\pi,q,h)\in \E_p(T)\,:\,h(0)=0\,\Big\}.\]
Restricted to this space the Fr{\'e}chet derivative of $N$ allows to be estimated as follows.
\begin{lemma}\label{lem:N}
	Let $\delta>0$, let $1 < p < \infty $ with $p \neq 4$, and let $T>0$.
	Then, $N\in C^\omega(C_\delta(T),\G_p(T))$, and for every fixed $z=(u,v,w,\pi,q,h)\in C_\delta(T)$ we have $DN(z)\in\mathcal{L}({_0\E_p}(T),{\G_p}(T))$ and 
	\begin{equation}\label{nichtlab}
		\begin{aligned}
			\|DN(z)\|_{\mathcal{L}({_0\E_p}(T),{\G_p}(T))}\le c\big(\|z\|_{\E_p(T)}+\|h\|_{C(\bar I,W^{5-4/p}_p(\Gamma))}\big),
		\end{aligned}
	\end{equation}
	where the constant $c>0$ is independent of $T$, but may depend on some upper bound for $\delta$.
\end{lemma}
\proof From \eqref{eqn:tildee} we see that pointwise all components of $N$ are analytic functions of $(u,v,w,\pi,q,h)$ and its derivatives. For the analyticity it thus suffices to prove that each term in $N:C_\delta(T)\rightarrow\G_p(T)$ is well-defined. This, however, is a rather simple exercise using Lemma \ref{lemma:emb}; cf. the proof of Proposition 6.2 in \cite{pruess62}. We present the idea by analyzing the most complicated nonlinearity $N_3^\perp$, leaving the other terms to the reader. Note that for dimensional reasons the terms in $Q(h)$ containing fourth order derivatives of $h$ must be of the form
\[r(h/\gamma,hk,\nabla h)*\nabla^4h.\]
By Lemma \ref{lemma:emb} (i) we have $\nabla h\in C(\bar I,C^2(\Gamma))$ which is, of course, an algebra with respect to pointwise multiplication. Since furthermoe $(\nabla^g)^4 h\in L_p(I,W^{1-1/p}_p(\Gamma))$ and
\begin{equation*}
\begin{aligned}
C(\bar I,C^2(\Gamma))\cdot L_p(I,W^{1-1/p}_p(\Gamma))\hookrightarrow L_p(I,W^{1-1/p}_p(\Gamma)),
\end{aligned}
\end{equation*}
that is, pointwise multiplication is continuous in the indicated function spaces, the terms containing fourth order derivatives of $h$ are well-defined. The terms involving $\nabla^g\tilde v$ contain up to second order derivatives of $h$. Since $\nabla^2 h\in C(I,C^1(\Gamma))$ which is also an algebra, $\nabla^g v\in L_p(I,H^{1}_p(\Gamma))$, and
\begin{equation*}
\begin{aligned}
C(\bar I,C^1(\Gamma))\cdot L_p(I,H^{1}_p(\Gamma))\hookrightarrow L_p(I,W^{1-1/p}_p(\Gamma)),
\end{aligned}
\end{equation*}
these terms are well-defined as well. The terms involving $q$ and $[\tilde u]$ can be handled analogously. Concerning the remaining terms in $Q(h)$ which contain up to third order derivatives of $h$ we simply note that, by Lemma \ref{lemma:emb} (i), $\nabla^3 h\in C(\bar I,W^{1-1/p}_p(\Gamma))$ and $C(\bar I,W^{1-1/p}_p(\Gamma))$ is an algebra for $p>3$; the latter fact follows from the theorem in Section 2.8.3 of \cite{Triebel:Function-Spaces-1} and a localization argument; cf. \cite{lengeler}. This completes the proof of analyticity for $N_3^\perp$. The other non-linearities can be handled analogously.

The estimate \eqref{nichtlab} essentially follows from the fact that $N$ vanishes with at least quadratic order in $z=0$; recall, in particular, the definition of $Q(h)$. The proof is again a rather simple exercise using Lemma \ref{lemma:emb} (cf. the proof of Proposition~4.1 in \cite{pruess41}) and, again, we present the idea by analyzing $DN_3^\perp$, leaving the other terms to the reader. All estimates derived below will be uniform in $T$. For some fixed $z=(u,v,w,\pi,q,h)\in C_\delta(T)$ and $\bar z=(\bar u,\bar v,\bar w,\bar\pi,\bar q,\bar h)\in {_0\E_p(T)}$ we have
\begin{equation*}
\begin{aligned}
DN_3^\perp(z)(\bar z)&=\tilde r(h/\gamma,hk,\nabla h)*(\nabla^g)^4 \bar h + \tilde r(h/\gamma,hk,\nabla h)*(\bar h/\gamma,\bar h k,\nabla\bar h)*(\nabla^g)^4 h\\ 
&\quad + \mbox{terms depending on up to third order derivatives of $h$ and $\bar h$}\\
&\quad +D\big(r(\tilde e)*\big((\tilde e-e)*k,\nabla\tilde e\big)\big)(h)(\bar h)\, q + r(\tilde e)*\big((\tilde e-e)*k,\nabla\tilde e\big)\,\bar q\\
&\quad +D\big(\mu\,r(\tilde e)*\big((\tilde e-e)*k,\nabla\tilde e\big)\big)(h)(\bar h)*\nabla^g v \\
&\quad + \mu\,r(\tilde e)*\big((\tilde e-e)*k,\nabla\tilde e\big)*\nabla^g\bar v\\
&\quad + D\big(\mu\,r(\tilde e)*\big((\tilde e-e)*k^2,k*\nabla\tilde e,(\nabla\tilde e)^2\big)\big)(h)(\bar h)*[u]\\
&\quad+ \mu\,r(\tilde e)*\big((\tilde e-e)*k^2,k*\nabla\tilde e,(\nabla\tilde e)^2\big)*[\bar u]
\end{aligned}
\end{equation*}
with analytic funtions $\tilde r$ such that $\tilde r(0,0,0)=0$. By the arguments used in the proof of analyticity, we have
\[\|r(h/\gamma,hk,\nabla h)*(\nabla^g)^4 \bar h\|_{L_p(I,W^{1-1/p}_p(\Gamma))}\le c\|h\|_{C(\bar I,C^2(\Gamma))}\,\|\bar z\|_{\E_p(T)}.\]
Similarly, using Lemma \ref{lemma:emb} (ii), we have
\begin{equation*}
\begin{aligned}
\|r(h/\gamma,hk,\nabla h)*(\bar h/\gamma,\bar h k,\nabla\bar h)*&(\nabla^g)^4 h\|_{L_p(I,W^{1-1/p}_p(\Gamma))}\\
&\le c\big(\|z\|_{\E_p(T)}+\|h\|_{C(\bar I,C^2(\Gamma))}\big)\|\bar h\|_{C(\bar I,C^2(\Gamma))}\\
&\le c\big(\|z\|_{\E_p(T)}+\|h\|_{C(\bar I,W^{5-4/p}_p(\Gamma))}\big)\|\bar z\|_{{_0\E}_p(T)}.
\end{aligned}
\end{equation*}
Again by the arguments used in the proof of analyticity, we have
\begin{equation*}
\begin{aligned}
\|r(\tilde e)*\big((\tilde e-e)*k,\nabla\tilde e\big)*\nabla^g\bar v\|_{L_p(I,W^{1-1/p}_p(\Gamma))}&\le c\|h\|_{C(J,C^3(\Gamma))}\,\|\bar z\|_{\E_p(T)}\\
&\le c\|h\|_{C(J,W^{5-4/p}_p(\Gamma))}\,\|\bar z\|_{\E_p(T)}
\end{aligned}
\end{equation*}
and, using Lemma \ref{lemma:emb} (ii),
\begin{equation*}
\begin{aligned}
\|D\big(\mu\,r(\tilde e)*\big((\tilde e-e)*k,\nabla\tilde e\big)\big)(h)(\bar h)*&\nabla^g v\|_{L_p(I,W^{1-1/p}_p(\Gamma))}\\
&\le c\big(\|z\|_{\E_p(T)}+\|h\|_{C(\bar I,C^3(\Gamma))}\big)\|\bar h\|_{C(I,C^3(\Gamma))}\\
&\le c\big(\|z\|_{\E_p(T)}+\|h\|_{C(\bar I,W^{5-4/p}_p(\Gamma))}\big)\|\bar z\|_{{_0\E}_p(T)}.
\end{aligned}
\end{equation*}
The terms involving $q$ and $[u]$ can be handled analogously. Finally, using
\[C(\bar I,W^{1-1/p}_p(\Gamma))\cdot L_p(I,W^{1-1/p}_p(\Gamma))\hookrightarrow L_p(I,W^{1-1/p}_p(\Gamma)),\]
we can estimate the terms depending only on up to third order derivatives of $h$ and $\bar h$ via
\begin{equation*}
\begin{aligned}
c\|h\|_{C(\bar I,W^{4-1/p}_p(\Gamma))}\|\bar h\|_{L_p(I,W^{4-1/p}_p(\Gamma))}\le c\|h\|_{C(\bar I,W^{5-4/p}_p(\Gamma))}\,\|\bar z\|_{\E_p(T)}.
\end{aligned}
\end{equation*}
This concludes the estimate of $DN^\perp_3(z)$. The derivatives of the other non-linearities can handled analogously.
\qed
\medskip

\proof[Proof of Theorem \ref{thm}] Following the remark after Definition \ref{def} we can show that $\G_p(T) = \F_p(T) \oplus \U_p(T)$ with
\[\U_p(T):=\Big\{\,(0,f_2,0,f_4,0,0,0)\in \G_p(T)\,:\,(f_2,f_4)\in L_p(I,U_p(\Gamma))\,\Big\}.\]
Let $P:\G_p(T)\rightarrow \F_p(T)$ denote the bounded projection along $\U_p(T)$. Furthermore, we write $L^{-1}:\F_p(T)\rightarrow{_0\E}_p(T)$ for the linear solution operator with $h(0)=0$ whose existence is guaranteed by Theorem \ref{thm:linear}. Since extension by $0$ defines a continuous operator $\F_p(T)\rightarrow \F_p(1)$ for $T<1$, we have a uniform bound
\[\|L^{-1}P\|_{\mathcal{L}(\G_p(T),{_0\E}_p(T))}\le M\]
for all $0<T\le 1$ and some $M>0$. From \eqref{nichtlab} and the inequality
\[\|h\|_{C(\bar I,W^{4-1/p}_p(\Gamma))}\le c\big(\|z\|_{\E_p(T)}+\|h_0\|_{W^{5-4/p}_p(\Gamma)}\big)\]
with a constant $c$ independent of $T$, by choosing $\delta$ and $\epsilon$ sufficiently small we obtain the estimate
\begin{equation}\label{DN}
\|DN(z)\|_{\mathcal{L}({_0\E}_p(T),\F_p(T))}\le \frac{1}{2M}
\end{equation}
for all $z\in C_\delta(T)$. Let $z^*=(u^*,v^*,w^*,\pi^*,q^*,h^*)\in \E_p(T)$ be the solution of $Lz^*=0$, $h^*(0)=h_0$ which exists according to Theorem \ref{thm:linear}; there exists a constant $c>0$ depending only on an upper bound for $T$ such that
\[\|z^*\|_{\E_p(T)}\le c\|h_0\|_{W^{5-4/p}_p(\Gamma)}.\]
We choose $\epsilon$ so small that $z^*\in C_{\delta/2}(T)$. Hence, we can write the transformed problem \eqref{eqn:tildesystem2} in the form
\begin{equation*}
 z=L^{-1}PN(z+z^*)=:K(z)
\end{equation*}
for some $z\in C'_{\delta/2}(T):=C_{\delta/2}(T)\cap{_0\E}_p(T)$. Note that $N(z^*)$ depends on $\grad_{L_2}F_\Gamma$ and $z^*$. Thus, in order to have $K(0)\in C'_{\delta/4}(T)$, we choose both $T$ and $\epsilon$, and hence $z^*\in \E_p(T)$, sufficiently small; the former choice has the effect that $\grad_{L_2}F_\Gamma$ is small in $L_p(I,W^{1-1/p}_p(\Gamma))$. By the contraction mapping principle the operator $K$ possesses a unique fixed point $z_0$ in $C'_{\delta/2}(T)$ if it maps this set contractively into itself. But this now follows from \eqref{DN} since we can infer
\begin{equation*}
\|DK(z)\|_{\mathcal{L}({_0\E}_p(T))}\le \frac{1}{2}
\end{equation*}
and
\[\|K(z)\|_{\E_p(T)}\le \|K(0)\|_{\E_p(T)} +  \frac12\|z\|_{\E_p(T)}\le \frac\delta 2\]
for all $z\in C'_{\delta/2}(T)$. Thus, for $\tilde z=(\tilde u,\tilde v,\tilde w,\tilde \pi,\tilde q,h)=z_0+z^*$ we have
\[L\tilde z=PN(\tilde z)=N(\tilde z)+(P-I)N(\tilde z);\]
note that $(P-I)N(\tilde z)=(0,f_2,0,f_4,0,0,0)$ for piecewise constant functions $f_2$ and $f_4$. Recalling the computations in \Secref{hanzawa} we see that the pushforward $u:=(\Phi_t)_*\tilde u$, $\pi:=(\Phi_t)_*\tilde \pi$, and $q:=(\Phi_t)_*\tilde q$ solves the system
\begin{equation}\label{phisystem}
\begin{aligned}
\Div S&=0 &&\mbox{ in }\Omega\setminus \Gamma_t,\\
\Div u&=(\Phi_t)_*f_2&&\mbox{ in }\Omega\setminus \Gamma_t,\\
\DIV {^fT}+[\![S]\!]\nu_t&=-\DIV {^eT}&&\mbox{ on }\Gamma_t,\\
\DIV u&=(\Phi_t)_*f_4&&\mbox{ on }\Gamma_t,\\
u&=0&&\mbox{ on }\pa\Omega
\end{aligned}
\end{equation}
for almost all $t\in I$, where the stress tensors are defined with respect to $u$, $\pi$, and $q$ and $\Gamma_t:=\Gamma_{h(t)}$. At this point we need to assume that $\Gamma$ contains no round spheres. Then, by definition of $U_p(\Gamma)$, we have $f_4=0$ and $f_2=\mbox{const}$ in $\Omega$. Now,  \eqref{phisystem}$_2$ shows that in fact $f_2=0$.

So far we proved that \eqref{eqn:final} has a local-in-time solution which is uniquely determined in the class of solutions whose transformation is of the form \[\tilde z=(\tilde u,\tilde v,\tilde w,\tilde \pi,\tilde q,h)=z_0+z^*\] with $z_0\in C'_{\delta/2}(T)$. Note that the pushforward and the pullback will in general not preserve the mean value condition \eqref{comp1} with $f_2=\tilde \pi(t,\cdot\,)$ on the one hand and with $f_2=\pi(t,\cdot\,)$ and $\Gamma_t$ in place of $\Gamma$ on the other hand; these conditions can be met, however, by adding suitable constants. Now, let us prove \emph{unconditional uniqueness} by a standard bootstrap argument; cf. \cite{Tao06}. To this end, let us repeat the contraction argument with $\delta/2$ in place of $\delta$ (leading to possibly smaller $\epsilon$ and $T$). We infer that in fact $z_0\in C'_{\delta/4}(T)$, but uniqueness still holds in $C'_{\delta/2}(T)$. Now, let $\tilde z'=(\tilde u',\tilde v',\tilde w',\tilde \pi',\tilde q',h')=z_0'+z^*$ with $z_0'\in {_0\E_
p(T)}$ denote another solution of \eqref{eqn:tildesystem2} with $\|h'\|_{L^\infty((0,T)\times \Gamma)}\le\kappa/2$; without restriction we may assume that it is defined on the same time interval as $\tilde z$. Choosing $T'$ sufficiently small we have $z_0'\in C'_{\delta/2}(T')$. Repeating again the contraction mapping argument, this time with $T'$ in place of 
$T$, we see that $z_0'$ coincides with $z_0$ on $(0,T')$; in particular we have $z_0'\in C'_{\delta/4}(T')$. But then we have $z_0'\in C'_{\delta/2}(T'')$ for some $T''$ slightly larger than $T'$. Thus, the set of times $T'$ with $z_0'\in C_{\delta/2}'(T')$ is open. But obviously it is also closed and non-empty, so that $z_0'\in C_{\delta/2}'(T)$ and $z_0'=z_0$. This proves unconditional uniqueness of our solution. 

Now, the proof of Lemma \ref{lem:N} shows that for fixed $T>0$  and all $z\in C_\delta(T)$ the norm $\|DN(z)\|_{\mathcal{L}(\E_p(T),\F_p(T))}$ is uniformly bounded. Hence, $N$ is Lipschitz continuous in $C_\delta(T)$, and thus for the operator $K=K_{h_0}$ we have
\[\|K_{h_0}(z)-K_{h_0'}(z)\|_{\E_p(T)}\le L\|K_{h_0}(z)-K_{h_0'}(z)\|_{W^{5-4/p}_p(\Gamma)}\]
for all $h_0,h_0'\in\bar B_\epsilon(0)\subset W^{5-4/p}_p(\Gamma)$, all $z\in C_{\delta/2}'(T)$, and some constant $L>0$. Now, for such initial values $h_0,h_0'$ let $z_{h_0}, z_{h_0'}\in C_{\delta/2}'(T)$ denote the respective fixpoints of $K_{h_0}$ and $K_{h_0'}$. Then, we have
\begin{equation*}
\begin{aligned}
\|z_{h_0} - z_{h_0'}\|_{\E_p(T)}&= \|K_{h_0}(z_{h_0}) - K_{h_0'}(z_{h_0'})\|_{\E_p(T)} \\ &\le  L \|h_0 - h_0'\|_{W^{5-4/p}_p(\Gamma)} + \frac12\|z_{h_0} -z_{h_0'}\|_{\E_p(T)}.
\end{aligned}
\end{equation*}
Absorbing the second term on the right hand side, we obtain the Lipschitz continuity of the solution map.

Finally, let us consider the case of $\Gamma$ being a collection of round spheres. Since the energy cannot decrease in this case, from \eqref{enid} we can show that $u$ must vanish everywhere, $\pi$ and $q$ are constant in each connected component of $\Omega$ and $\Gamma$, respectively, and
\begin{equation*}
 \kappa\frac{C_0}{R_i}\bigg(\frac{2}{R_i} -C_0\bigg)+[\![\pi]\!]+q\,\frac{2}{R_i}=\grad_{L_2}F+[\![\pi]\!]+q\,H=0
\end{equation*}
on $\Gamma^i$, $i=1,\ldots,m$, where $\Gamma^i$ is a round sphere of radius $R_i$; for details see the discussion in the end of Section 2 of \cite{lengeler}. Combining these $m$ conditions with the $m+1$ conditions \eqref{comp1} and \eqref{comp2} we obtain a system of linear equations which can easily be uniquely solved for the $2m+1$ unkowns $q$ on $\Gamma^i$, $\pi$ in $\Omega^i$, and $\pi$ in $\Omega^0$.
\qed
\medskip

Suppose that each $\Gamma^i$, $i=1,\ldots,l$ and $l\ge 1$, is a round sphere while each $\Gamma^i$, $i=l+1,\ldots,m$ and $m\ge 2$, is a non-sphere and that $h_0=0$. Then, in general, the potential solution will not be constant in time and the round spheres might translate. In this case, however, showing that $f_2$ and $f_4$ in the above proof vanish is not completely obvious. If we know that at some fixed instant $t$ in time each $\Gamma_t^i$, $i=1,\ldots,l$, is a round sphere, then by definition of $U_p(T)$ it is not hard to see that $f_2$ and $f_4$ must vanish at time $t$. Thus, by \eqref{phisystem}$_{2,4}$, each $\Gamma_t^i$, $i=1,\ldots,l$, will remain a round sphere for the next instant in time (in linear approximation). This situation suggests to apply some kind of continuity or Gronwall-type argument; so far, however, we were not able to close the required estimates. On the other hand, it's questionable if this slight generalization of our theorem is worth the effort.

\appendix
\section{Covariant differentiation and curvature}\applabel{app:cov}
Here, we recall some useful results from Appendix~B in \cite{lengeler}. Let $e_{ij}$, $\tilde e_{ij}$ be Riemannian metrics on a manifold $M$, and let $e^{ij}$, $\tilde e^{ij}$ denote their matrix inverses. For scalar functions $f$, vector fields $Y$, and second order tensor fields $T$ we have
\begin{equation*}
\begin{aligned}
(\grad_{\tilde e}f)^i&=(\grad_{e}f)^i + (\tilde e^{ij}-e^{ij})\pa_j f,\\
\Div_{\tilde e}Y&=\Div_e Y + \tilde e*\nabla^e\tilde e*Y,\\
\Delta_{\tilde e} f &= \Delta_e f + (\tilde e-e)*r(\tilde e,e)*(\nabla^e)^2 f + r(\tilde e,e)*\nabla^e\tilde e*\nabla f,\\
D^{\tilde e}Y&=D^{e}Y+(\tilde e-e)*r(\tilde e,e)*\nabla^e Y + r(\tilde e,e)*\nabla^e \tilde e *Y,\\
\Delta_{\tilde e} Y&=\Delta_{e} Y + (\tilde e-e)*r(\tilde e,e)*(\nabla^e)^2 Y + r(\tilde e,e)*(\nabla^e)^2\tilde e*Y\\ 
&\quad+ r(\tilde e,e)*(\nabla^e\tilde e)^2*Y
+ r(\tilde e,e)*\nabla^e\tilde e*\nabla^e Y,\\
\Div_{\tilde e} T&=\Div_{e} T + (\tilde e-e)*r(\tilde e,e)*\nabla^e T + r(\tilde e,e)*\nabla^e\tilde e*T.
\end{aligned}
\end{equation*}
where $D^e Y$ and $D^{\tilde e} Y$ denote the $e$-symmetric part of $\nabla^e Y$ and the $\tilde e$-symmetric part of $\nabla^{\tilde e}Y$, respectively. Furthermore, let $\Gamma$ be an orientable submanifold of $M$ of codimension $1$, and let $\nu_{e}$ and $\nu_{\tilde e}$ be equally oriented unit normal fields on $\Gamma$ with respect to $e$ and $\tilde e$, respectively. Then, we have
\begin{equation*}
\begin{aligned}
\nu_{\tilde e}&=\nu_e + (\tilde e-e)*r(\tilde e,e),\\
\nabla^{\tilde e} \nu_{\tilde e}&= \nabla^e\nu_{e} + r(\tilde e,e)*\nabla^e\tilde e,\\
k_{\tilde e}&=k_e+(\tilde e-e)*k_e + r(\tilde e,e)*\nabla^e\tilde e,\\
H_{\tilde e}&= H_e + (\tilde e-e)*r(\tilde e,e)*k_e + r(\tilde e,e)*\nabla^e\tilde e,\\
K_{\tilde g}&=\det(\tilde g^{\alpha\delta}(k_{\tilde e})_{\delta\beta})\\
&=\det\big(g^{\alpha\delta}(k_e)_{\delta\beta} + (\tilde e-e)*r(\tilde e,e)*k_e + r(\tilde e,e)*\nabla^e\tilde e\big)\\
&=K_g+r(\tilde e,e)*\big((\tilde e-e)*k^2_e,k_e*\nabla\tilde e,(\nabla\tilde e)^2\big).
\end{aligned}
\end{equation*}

\section{The Stokes system in $\bR^n$ and $\bR^n_+$}
\seclabel{appendix:stokes}
Let $1 < p < \infty$. We consider the stationary Stokes system
\begin{equation*}
	\begin{array}{rcll}
		\eta\, u - \mu \Delta u + \grad \pi & = & f & \quad \mbox{in} \ \bR^n, \\[0.5em]
		                   \mbox{div}\,u & = & g & \quad \mbox{in} \ \bR^n
	\end{array}
\end{equation*}
for some shift $\eta > 0$ and some constant viscosity $\mu > 0$. There exists a unique solution
\begin{equation*}
	u \in H^2_p(\bR^n,\bR^n), \quad \pi \in \dot{H}^1_p(\bR^n)/\R,
\end{equation*}
provided that $f \in L_p(\bR^n,\,\bR^n)$, and $g \in H^1_p(\bR^n)$. Indeed, we may first obtain the pressure as $\pi = (-\mu + \eta (-\Delta)^{-1}) g - \mbox{div}\,(-\Delta)^{-1} f \in \dot{H}^1_p(\bR^n)$ to be left with the equation
\begin{equation*}
	\eta\, u - \mu \Delta u = f - \grad p \quad \mbox{in} \ \bR^n,
\end{equation*}
which allows for a solution $u \in H^2_p(\bR^n,\bR^n)$, since this is an elliptic problem with right-hand-side $f - \grad\pi \in L_p(\bR^n,\bR^n)$.
Finally, uniqueness of solutions is a direct consequence of the validity of the Helmholtz decomposition in $L_p(\R^n,\R^n)$.

As a direct consequence, we infer that the Stokes system
\begin{equation*}
	\begin{array}{rcll}
		\eta\, u - \mu \Delta u + \grad \pi & = & f      & \quad \mbox{in} \ \bR^n_+, \\[0.5em]
		                   \mbox{div}\,u & = & g_p    & \quad \mbox{in} \ \bR^n_+, \\[0.5em]
		                      [v]_\Sigma & = & g_\tau & \quad \mbox{on} \ \Sigma,  \\[0.5em]
		                      [\pi]_\Sigma & = & g_\nu  & \quad \mbox{on} \ \Sigma
	\end{array}
\end{equation*}
in the halfspace $\bR^n_+ := \{\,(x,\,y) \in \bR^{n - 1} \times \bR\,:\,y > 0\,\}$ also allows for a unique solution
\begin{equation*}
	u \in H^2_p(\bR^n_+,\bR^n), \quad \pi \in \dot{H}^1_p(\bR^n_+)/\R,
\end{equation*}
provided that $f \in L_p(\bR^n_+,\bR^n)$, $g_p \in H^1_p(\bR^n_+)$, $g_\tau \in W^{2 - 1/p}_p(\Sigma,\,\bR^{n - 1})$, as well as $g_\nu \in \dot{W}^{1 - 1/p}_p(\Sigma)$.
We employ the usual decomposition \mbox{$u = (v,\,w)\in\R^{n-1}\times\R$} and note that the trace operator
\begin{equation*}
	[\,\cdot\,]_\Sigma: \dot{H}^1_p(\bR^n_+) \longrightarrow \dot{W}^{1 - 1/p}_p(\Sigma)
\end{equation*}
admits a bounded linear right-inverse as follows from \cite[Theorems 2.4' and 2.7, Corollary 1]{Kudrjavcev:Imbedding-1, Kudrjavcev:Imbedding-2}.
Now, we may first eliminate $g_\tau$ and $g_\nu$ by constructing extensions $\bar{v} \in H^2_p(\bR^n_+,\,\bR^{n - 1})$ to $g_\tau$ and $\bar{\pi} \in \dot{H}^1_p(\bR^n_+)$ to $g_\nu$ and then solve the remaining problem by a reflection argument via a Stokes problem in $\bR^n$; more precisely, for $f=(f_1,\ldots,f_n)$, we extend $f_1,\ldots,f_{n-1}$ and $g_p$ by an odd reflection and $f_n$ by an even reflection to $\R^n$ .

As another consequence, we infer that the Stokes system
\begin{equation*}
	\begin{array}{rcll}
		\eta\, u - \mu \Delta u + \grad\pi & = & f      & \quad \mbox{in} \ \bR^n_+, \\[0.5em]
		                   \mbox{div}\,u & = & g_p    & \quad \mbox{in} \ \bR^n_+, \\[0.5em]
		                      [v]_\Sigma & = & g_\tau & \quad \mbox{on} \ \Sigma,  \\[0.5em]
		                      [w]_\Sigma & = & g_\nu  & \quad \mbox{on} \ \Sigma
	\end{array}
\end{equation*}
allows for a unique solution
\begin{equation*}
	u \in H^2_p(\bR^n_+,\bR^n), \quad \pi \in \dot{H}^1_p(\bR^n_+)/\R,
\end{equation*}
too, provided that $f \in L_p(\bR^n_+,\bR^n)$, $g_p \in H^1_p(\bR^n_+)$, $g_\tau \in W^{2 - 1/p}_p(\Sigma,\bR^{n - 1})$, and $g_\nu \in W^{2 - 1/p}_p(\Sigma)$. Indeed, we may in a first step eliminate $f$ and $g_p$ by extending these function to $\bR^n$ and solving the corresponding Stokes system in the whole space. The reduced problem may then be treated with the aid of a Fourier transform in the tangential variables $x \in \bR^{n - 1}$, that is, we consider the system
\begin{equation*}
	\begin{array}{rcll}
		     \eta\, \hat{v} + \mu |\xi|^2 \hat{v} - \mu \partial^2_y \hat{v} + i \xi \hat{\pi} & = & 0            & \quad \xi \in \bR^{n - 1}, \ y > 0, \\[0.5em]
		\eta\, \hat{w} + \mu |\xi|^2 \hat{w} - \mu \partial^2_y \hat{w} + \partial_y \hat{\pi} & = & 0            & \quad \xi \in \bR^{n - 1}, \ y > 0, \\[0.5em]
		                                             i \xi^{\sfT} \hat{v} + \partial_y \hat{w} & = & 0            & \quad \xi \in \bR^{n - 1}, \ y > 0, \\[0.5em]
		                                                                      [\hat{v}]_\Sigma & = & \hat{g}_\tau & \quad \xi \in \bR^{n - 1}, \ y = 0, \\[0.5em]
		                                                                      [\hat{w}]_\Sigma & = & \hat{g}_\nu  & \quad \xi \in \bR^{n - 1}, \ y =0 ,
	\end{array}
\end{equation*}
The solution again has the form \eqref{schwuff}, and a straight forward computation yields
\begin{equation*}
	\left[ \begin{array}{c} \hat{z}_v(\xi) \\[0.5em] \hat{z}_w(\xi) \end{array} \right]
		= \frac{1}{\varpi} \left( \left(1 - \frac{|\zeta|}{\varpi} \right) \frac{|\zeta|}{\varpi} \right)^{-1}
			\left[ \begin{array}{rr} \left(1 - \frac{|\zeta|}{\varpi} \right) \frac{|\zeta|}{\varpi} - \frac{i \zeta \otimes i \zeta}{\varpi^2} & \frac{i \zeta}{\varpi} \\[0.5em] \quad - \frac{i \zeta^{\sfT}}{\varpi} & \quad 1 \end{array} \right]
			\left[ \begin{array}{c} \hat{g}_\tau(\xi) \\[0.5em] \hat{g}_\nu(\xi) \end{array} \right];
\end{equation*}
in particular, we have
\begin{equation*}
	\begin{aligned}
		\widehat{\grad_x \pi}(\xi,\,y)
			& =  \eta \sqrt{\mu_b}\,i \xi\,\hat{z}_w(\xi) e^{- |\xi| y} \\[1.5em]
			& = \eta\,\frac{i \zeta}{\varpi} \left( \left(1 - \frac{|\zeta|}{\varpi} \right) \frac{|\zeta|}{\varpi} \right)^{-1} \left( \hat{g}_\nu - \frac{i \zeta^{\sfT}}{\varpi} \hat{g}_\tau \right) e^{- |\xi| y} \\[1.5em]
			& = \sqrt{\mu}\,\frac{\varpi}{|\zeta|} \frac{i \zeta}{|\zeta|} (\varpi + |\zeta|) \left( \hat{g}_\nu - \frac{i \zeta^{\sfT}}{\varpi} \hat{g}_\tau \right) |\xi| e^{- |\xi| y}
	\end{aligned}
\end{equation*}
and
\begin{equation*}
	\widehat{\partial_y \pi}(\xi,\,y)
		= - \eta \sqrt{\mu}\,|\xi|\,\hat{z}_w(\xi) e^{- |\xi| y}
		= - \sqrt{\mu}\,\frac{\varpi}{|\zeta|} (\varpi + |\zeta|) \left( \hat{g}_\nu - \frac{i \zeta^{\sfT}}{\varpi} \hat{g}_\tau \right) |\xi| e^{- |\xi| y}.
\end{equation*}
Here, the symbol $|\xi| e^{- |\xi| y}$ belongs to the operator $A T(y)$,
where $A = (-\Delta)^{1/2}$ and $T(\,\cdot\,)$ denotes the corresponding semigroup.
Since $- \Delta_x q = A T(\,\cdot\,) h$ for the unique solution $q \in \dot{H}^2_p(\bR^n_+)$ of the elliptic boundary value problem
\begin{equation*}
\begin{aligned}
	- \Delta q &= 0 \quad \mbox{in}\ \bR^n_+,\\
	\partial_\nu q &= h \quad \mbox{on}\ \Sigma
\end{aligned}
\end{equation*}
with $h \in \dot{W}^{1 - 1/p}_p(\Sigma)$, we infer that
\begin{equation*}
	A T(\,\cdot\,): \dot{W}^{1 - 1/p}_p(\Sigma) \rightarrow L_p(\bR^n_+)
\end{equation*}
is bounded. Combining this observation with Mikhlin's multiplier theorem we conclude that $\pi \in \dot{H}^1_p(\bR^n_+)$. Then, the velocity field may be obtained as a solution of the elliptic boundary value problem
\begin{equation*}
	\begin{array}{rcll}
		\eta\, u - \mu \Delta u & = & f - \grad \pi  & \quad \mbox{in} \ \bR^n_+, \\[0.5em]
		             [v]_\Sigma & = & g_\tau         & \quad \mbox{on} \ \Sigma,  \\[0.5em]
		           [w]_\Sigma & = & g_\nu            & \quad \mbox{on} \ \Sigma,
	\end{array}
\end{equation*}
which implies $u \in H^2_p(\bR^n_+,\,\bR^n)$.


\bibliographystyle{plain}
\bibliography{references}
\end{document}